\documentclass[12pt]{article}
\usepackage{graphicx}
\usepackage{amsmath,amsthm,amssymb,enumerate}
\usepackage{euscript,mathrsfs}
\usepackage{amsmath}
\usepackage{yfonts}
\usepackage{color}
\usepackage{dsfont}
\usepackage{mathrsfs}
\usepackage[left=2cm,right=2cm,top=3.5cm,bottom=3.5cm]{geometry}
\usepackage[framemethod=tikz]{mdframed}
\usepackage{bm}
\allowdisplaybreaks

\usepackage{comment}

\usepackage{soul}

\catcode`\@=11 \@addtoreset{equation}{section}

\catcode`\@=12

\allowdisplaybreaks

\newtheorem{Theorem}{Theorem}[section]
\newtheorem{Proposition}[Theorem]{Proposition}
\newtheorem{Lemma}[Theorem]{Lemma}
\newtheorem{Corollary}[Theorem]{Corollary}

\newtheorem{Definition}[Theorem]{Definition}

\newtheorem{Remark}[Theorem]{Remark}

\newcommand{\bTheorem}[1]{
\begin{Theorem} \label{T#1} }
\newcommand{\eT}{\end{Theorem}}

\newcommand{\bProposition}[1]{
\begin{Proposition} \label{P#1}}
\newcommand{\eP}{\end{Proposition}}

\newcommand{\bLemma}[1]{
\begin{Lemma} \label{L#1} }
\newcommand{\eL}{\end{Lemma}}

\newcommand{\bCorollary}[1]{
\begin{Corollary} \label{C#1} }
\newcommand{\eC}{\end{Corollary}}

\newcommand{\bRemark}[1]{
\begin{Remark} \label{R#1} }
\newcommand{\eR}{\end{Remark}}

\newcommand{\bDefinition}[1]{
\begin{Definition} \label{D#1} }
\newcommand{\eD}{\end{Definition}}

\newcommand{\Del}{\Delta_x}

\newcommand{\Ds}{\mathbb{D}_x}

\newcommand{\CH}{\mathds{1}_{\mathcal{O}_{\rm ess}} }

\newcommand{\vuB}{\vc{u}_B}

\newcommand{\bfphi}{\boldsymbol{\varphi}}

\newcommand{\bFormula}[1]{
\begin{equation} \label{#1}}
\newcommand{\eF}{\end{equation}}

\newcommand{\vrn}{\vr_n}
\newcommand{\vun}{\vu_n}
\newcommand{\vtn}{\vt_n}

\newcommand{\Ov}[1]{\overline{#1}}

\newcommand{\aleq}{\stackrel{<}{\sim}}

\newcommand{\ageq}{\stackrel{>}{\sim}}

\newcommand{\vr}{\varrho}

\newcommand{\tvr}{\tilde \vr}
\newcommand{\tvu}{{\tilde \vu}}
\newcommand{\tvt}{\tilde \vt}

\newcommand{\vt}{\vartheta}
\newcommand{\vu}{\vc{u}}
\newcommand{\vm}{\vc{m}}

\newcommand{\vc}[1]{{\bm #1}}

\newcommand{\Div}{{\rm div}_x}
\newcommand{\Grad}{\nabla_x}

\newcommand{\dx}{\,{\rm d} {x}}

\newcommand{\dt}{\,{\rm d} t }

\newcommand{\intO}[1]{\int_{\Omega} #1 \ \dx}
\newcommand{\intOn}[1]{\int_{\Omega_n} #1 \ \dx}

\newcommand{\intR}[1]{\int_{\R^3} #1 \ \dx}

\newcommand{\D}{{\rm d}}

\newcommand{\R}{\mathbb{R}}
\newcommand{\vtB}{\vt_B}

\newcommand{\br}{ \nonumber \\ }

\newcommand{\bb}[1]{{\mathbb #1}}
\def\softd{{\leavevmode\setbox1=\hbox{d}%
	\hbox to 1.05\wd1{d\kern-0.4ex{\char039}\hss}}}
\definecolor{Cgrey}{rgb}{0.85,0.85,0.85}
\definecolor{Cblue}{rgb}{0.50,0.85,0.85}
\definecolor{Cred}{rgb}{1,0,0}
\definecolor{fancy}{rgb}{0.10,0.85,0.10}

\newcommand{\cblue}{\color{blue}}

\newcommand\Cbox[2]{%
\newbox\contentbox%
\newbox\bkgdbox%
\setbox\contentbox\hbox to \hsize{%
	\vtop{
		\kern\columnsep
		\hbox to \hsize{%
			\kern\columnsep%
			\advance\hsize by -2\columnsep%
			\setlength{\textwidth}{\hsize}%
			\vbox{
				\parskip=\baselineskip
				\parindent=0bp
				#2
			}%
			\kern\columnsep%
		}%
		\kern\columnsep%
	}%
}%
\setbox\bkgdbox\vbox{
	\color{#1}
	\hrule width  \wd\contentbox %
	height \ht\contentbox %
	depth  \dp\contentbox
	\color{black}
}%
\wd\bkgdbox=0bp%
\vbox{\hbox to \hsize{\box\bkgdbox\box\contentbox}}%
\vskip\baselineskip%
}

\mdfdefinestyle{MyFrame}{%
linecolor=black,
outerlinewidth=1pt,
roundcorner=5pt,
innertopmargin=\baselineskip,
innerbottommargin=\baselineskip,
innerrightmargin=10pt,
innerleftmargin=10pt,
backgroundcolor=white!20!white}

\date{}

\makeindex
\begin{document}


\title{On thermally driven fluid flows arising in astrophysics}

\author{Nilasis Chaudhuri\thanks{The research of NC was supported by the EPSRC Early Career Fellowship  EP/V000586/1, and by the ``Excellence Initiative Research University (IDUB)" program at the University of Warsaw.} \and Eduard Feireisl
	\thanks{The work of E.F. was supported by the
		Czech Sciences Foundation (GA\v CR), Grant Agreement
		24-11034S. The Institute of Mathematics of the Academy of Sciences of
		the Czech Republic is supported by RVO:67985840. } \and  Ewelina Zatorska\thanks{The research of  EZ  was supported by the EPSRC Early Career Fellowship  EP/V000586/1.} \and Bogus\l aw Zegarli\'nski$^{\mathsection}$ 
}

\date{\today}

\maketitle

\bigskip
\centerline{$^*$ {Faculty of Mathematics, Informatics, and Mechanics, University of Warsaw,}} \centerline{ ul. Banacha 2, Warsaw 02-097, Poland}

\centerline{$^\dag$   Institute of Mathematics of the Academy of Sciences of the Czech Republic}

\centerline{\v Zitn\' a 25, CZ-115 67 Praha 1, Czech Republic}

\centerline{$^\ddag$
Mathematics Institute, University of Warwick, Coventry CV4 7AL, United Kingdom}

\centerline{$^{\mathsection}$   
	Institute de Mathematiques de Toulouse, CNRS,   UMR5219 }
 
\begin{abstract}
	
	We consider the Navier-Stokes-Fourier-Poisson system  driven by an inhomogeneous temperature distribution on the boundary of an exterior fluid domain. We impose the finite mass constraint, positive far field condition for the temperature as well as  the no--slip boundary conditions for the velocity. The existence of global--in--time weak solutions and the weak-strong principle are proved.

\end{abstract}


\bigskip

\tableofcontents
%
%

\section{Introduction} 
\label{i}

We are concerned with the global-in-time existence theory for weak solutions of the multidimensional compressible Navier-Stokes-Fourier-Poisson (NSFP) system with large initial data. Specifically, we consider the NSFP system in the domain exterior to a bounded subset of $\R^3$, assuming that the mass of the fluid is finite.  Our physical motivation is to propose a possible description of the motion of compressible gaseous stars under a self-consistent gravitational field. 
Understanding such a system in its full complexity poses considerable challenges. Therefore, to gain insight and develop technical tools for further study, it is necessary to propose and to study simplifying physically realistic models that capture essential initial features.
For instance, when considering a star, one may initially overlook complications arising from the plethora of interesting and important phenomena in its upper non-solid part. Instead, as we propose, one can study a three-dimensional model with specified boundary conditions on the surface of the hard core of the star, to describe the flow in the unbounded outer region. This approach serves as a foundational model that can later be expanded to incorporate effects such as magnetic fields, rotation of the solid core, or substance inflow-outflow, as demonstrated in \cite{ChNY, FeGwKwSG} and related works.

The key point of our analysis is a proper choice of boundary/far field condition.
Specifically, we impose the zero far field condition for the density compatible with finite mass of the gas. 
The fluid motion is driven the prescribed temperature on the interface between the solid core and the gaseous atmosphere.
Additionally, our model allows for the solid core of the star to be of arbitrary shape. In contrast with other results assuming radial symmetry \cite{CWY, Jiang, Makino} or thermal insulation \cite{Poul1} of the outer domain, our setup is tailored to spherically inhomogeneous, thermally driven, and non-stationary astrophysical flows.

Our goal is to extend the mathematical theory of open fluid systems developed in 
\cite{ChauFei} and \cite{FeiNovOpen} to problems on \emph{unbounded} (exterior) domains, in particular to the NSFP system (formulated in equations \eqref{NSF1}-\eqref{FL} below). One of the principal difficulties is the finite mass constraint forcing the density to vanish in the far field. 
Accordingly, the momentum and energy equations may degenerate in the far field, loosing their ''parabolic'' character. 
Similarly to the classical Rayleigh--B\' enard problem, the fluid is pushed to the outer space by the buoyancy force and the balance can be restored 
only by gravitation. To the best of our knowledge, this type of problem has never been rigorously analysed in the literature.

Similarly to \cite{JeJiNo} and \cite{Poul}, we 
construct solutions by means of a sequence of approximate problems on bounded 
invading domains $(\Omega_n)_{n=1}^\infty$, specifically, 
\[
\Omega_n = \Omega \cap \left\{ x \in \R^3 \ \Big|\ |x| \leq r + n \right\},
\]
with suitable conditions imposed on the outer boundary. The existence of approximate solutions $(\vrn, \vtn, \vun)_{n=1}^\infty$ 
was proved in \cite{ChauFei} or \cite[Chapter 12]{FeiNovOpen}. 
As the compactness arguments necessary to perform the asymptotic limit $n \to \infty$ are essentially local and nowadays well understood (cf. \cite{JeJiNo}, \cite{Poul}), the main issue is to find suitable uniform bounds on the sequence of approximate solutions. These can be obtained by careful analysis of the associated ballistic energy functional. 

Once the existence of weak solutions is established,
we extend the scope of the relative energy method \cite{FeiNov2021NON} to prove the weak-strong uniqueness property. Specifically, the weak and strong solutions emanating from the same data coincide on the lifespan of the strong solution.
As a matter of fact, the issue is quite delicate as the existence of a (local) strong solutions for this particular model is nor known. Indeed the 
main difficulty is the combination of the general geometry of the fluid domain, the Dirichlet boundary conditions of the temperature, and the assumption of
finite mass. Some  partial results for small data solutions can be found in \cite{ChoKim1, WenZhu2}, and the references \cite{DanchinMucha, HuLiXi} address the isentropic case. Other existing results with similar far field conditions for density include  \cite{ChoKim1}, \cite{HuLiXi}, and \cite{WenZhu2}.
In \cite{LX, LX_2022A} it was shown that  the absolute temperature becomes uniformly positive at each positive time, no matter whether it is uniformly positive or not initially. This result was however restricted to the one dimensional compressible Navier-Stokes-Fourier equations  with density decaying at infinity  sufficiently fast. For density and temperature decaying at similar rates, the well-posedness  for the Cauchy problem in the multi-dimensional setting was  obtained in \cite{LX_2022B}. More recently, in \cite{DXZ}, analogous results were obtained for the system with degenerate viscosities and heat conductivity depending on the absolute temperature in a power law. They key ingredient was to rewrite the system in terms of entropy rather than temperature and to show its  uniformly finite lower and upper bounds. We anticipate the regularity of the local-in-time regular solution based on these results rewritten in terms of the absolute temperature $\vt$.

\subsection{Field equations}
In our description of fluid motion, we adhere to the framework and terminology of continuum fluid mechanics. Accordingly, the movement of a compressible, viscous, and heat-conducting fluid is described in terms of three basic phase variables: mass density $\vr = \vr(t,x)$, bulk velocity $\vu = \vu(t,x)$, and  (absolute) temperature $\vt= \vt(t,x)$. 
The time evolution of the fluid is governed by the Navier--Stokes--Fourier (NSF) system -- a mathematical formulation of 
the mass and  linear momentum conservation, and of the internal energy balance, respectively:
\begin{align}
	&\partial_t \varrho + \Div (\varrho \vc{u}) =0,\label{NSF1}\\
	&\partial_t (\varrho \vc{u}) + \Div (\varrho \vc{u} \otimes \vc{u}) + \nabla_x p(\varrho,\vartheta) = \Div \bb{S}(\vt , \Ds \vu) + \varrho  \vc{g}, \label{NSF2}\\
	&\partial_t (\varrho e(\varrho,\vartheta)) + \Div (\varrho e(\varrho,\vartheta)\vu ) + \Div \vc{q}(\vt, \Grad \vt) = \bb{S}(\vt, \Ds \vu) \colon \bb{D}_x \vc{u} - p(\varrho , \vartheta) \Div \vc{u} ,\label{NSF3}
\end{align}
where  the pressure $p = p(\vr, \vt)$ and the internal energy $e = e(\vr, \vt)$ are determined in terms $( \vr,\vt )$ through suitable equations of state (EOS). 

The symbol $\vc{g} = \vc{g}(t,x)$  represents the gravitation. Specifically, 
\begin{equation} \label{G1}
	\vc{g} =  \vr \Grad G,\quad - \Div (\mathfrak{g}\Grad G) = \vr + \Ov{\vr}\quad {\text{in}}\ \R^3,
\end{equation}
where $\mathfrak{g} > 0$ is the gravitational constant, $\vr$ is the density of the gas (auto--gravitation), 
and $\Ov{\vr}$ is the mass density of the ``hard core'' supported on 
$\R^3 \setminus \Omega$, see e.g. 
\cite{DFPS1}.

The viscous stress tensor $ \mathbb{S} $ is given by \emph{Newton's rheological law} 
	\begin{align}\label{NR}
	\mathbb{S}=\mathbb{S}(\vt,\Ds \vu)=2\mu( \vartheta) \bigg(\Ds \vc{u}-\frac{1}{3} (\Div\vc{u})\mathbb{I} \bigg) + \lambda( \vartheta) (\Div \vc{u}) \mathbb{I},
\end{align}
where $\ \Ds \vu \equiv \frac{1}{2} \left( \Grad \vu + \Grad^t \vu \right)$ is the symmetric velocity gradient,  $ \mu(\vt) > 0 $ is the \emph{shear viscosity} coefficient, and $ \lambda(\vt) \geq 0 $ is the \emph{bulk viscosity} coefficient.

The heat flux $\vc{q}$ is given by \emph{Fourier's law}, 
\begin{align}\label{FL}
	\vc{q}(\vt, \Grad \vt) = - \kappa (\vt) \nabla_x \vartheta,
\end{align}
with the \emph{heat conductivity} coefficient $ \kappa(\vt) > 0 $. 

\subsection{Second Law of Thermodynamics, entropy}

The EOS determining the pressure $p= p(\varrho,\vartheta) $ and the internal energy $ e=e(\varrho, \vartheta) $ comply with \emph{Gibbs' equation},
\begin{align}\label{Gib}
	\vartheta Ds= De+pD\left(\frac{1}{\varrho}\right),
\end{align}
where the quantity $ s=s(\varrho,\vartheta)  $ is called the \emph{entropy}. Accordingly, 
the \emph{internal energy balance} \eqref{NSF3} can be rewritten as the \emph{entropy equation} 
\begin{align}
		\partial_t (\vr s(\vr, \vt)) + \Div (\vr s(\vr, \vt) \vu) + \Div \left( \frac{ \vc{q}(\vt, \Grad \vt) }{\vt} \right) = \frac{1}{\vt} 
		\left( \mathbb{S}(\vt, \Ds \vu)  : \bb{D}_x \vu - \frac{\vc{q}(\vt , \Grad \vt) \cdot \nabla_x \vt }{\vt} \right). \label{ent_eq}
\end{align}
It follows from \eqref{NR}, \eqref{FL} that the entropy production rate is non-negative, 
\[
\frac{1}{\vt} 
\left( \mathbb{S}(\vt, \Ds \vu)  : \bb{D}_x \vu - \frac{\vc{q} (\vt, \Grad \vt) \cdot \nabla_x \vt }{\vt} \right) \geq 0,
\]
in agreement with the Second Law of Thermodynamics.

\subsection{Fluid domain and the boundary conditions}

Inspired by models in astrophysics, we suppose the fluid is confined to an exterior domain $\Omega \subset \R^3$; more specifically, 
$\Omega$ is an unbounded domain with a compact boundary of class $C^{2 + \nu}$.
The compact complement $\R^3 \setminus \Omega$ represents a hard core, where the fluid equations are not relevant.  

Our main working hypothesis asserts the total mass of the 
fluid 
\begin{equation} \label{BC1a}
m_0 = \intO{ \vr(t, \cdot) }.
\end{equation}
to be constant and finite.

The total mass of the fluid is conserved under the impermeability condition 
\[
\vu \cdot \vc{n}|_{\partial \Omega} = 0,  	
\]	
where $\vc{n}$ denotes the outer normal to $\partial \Omega$ hereafter.
However, to determine the time evolution of the tangential part of the velocity, complementary conditions are necessary. For the sake of simplicity, 
we consider the no--slip condition:  
\begin{equation} \label{BC1}
	\vc{u}|_{\partial \Omega} = \vc{0}.
	\end{equation}
\begin{Remark}
   We refer \cite{ChauFei} for the necessary modifications of the proof to accommodate the inhomogeneous 
   	boundary conditions for the velocity $\vu = \vuB$, $\vuB \cdot \vc{n} = 0$.
   
\end{Remark}	
In addition, we prescribe the boundary temperature, 
\begin{equation} \label{BC2}
	\vt|_{\partial \Omega} = \vtB,\quad \vtB=\vtB(t,x) > 0.
\end{equation}
	
Finally, we impose the far field conditions
	\begin{align}
	&\vr \rightarrow 0,\; \vu \rightarrow \vc{0} \text{ and } \vt \rightarrow \vt_\infty \text{ as } \vert x \vert \rightarrow \infty, 
		\br
	&\mbox{where}\ \vt_\infty \ \mbox{is constant},\ 0 < \vt_\infty \leq \inf_{\partial \Omega} \vtB.	
	\label{FFC2}
\end{align}

On the one hand, 
note that \eqref{FFC2} is compatible with the total mass of the fluid $m_0$ determined
through  \eqref{BC1a} to be finite. On the other hand, the temperature $\vt$ satisfies the inhomogenous and possible time-dependent boundary conditions 
on $\partial \Omega$ supplemented with a positive far field condition stated in \eqref{FFC2}. To the best of our knowledge, this combination of boundary conditions has never been considered in the literature, neither in the weak nor in the strong solutions framework.  There are small data global existence results 
in the case $\Omega = \R^3$, $\vt_\infty = 0$ by Cho and Kim 
\cite{ChoKim1}, and by Wen and Zhu \cite{WenZhu2}, among others. The existence of large data weak solutions was established by
Poul \cite{Poul1} in the case of thermally isolated exterior domain, and by 
Jessl\' e, Jin, and Novotn\' y \cite{JeJiNo} in the case of thermally isolated boundary and positive far field density 
\[
\vr \to \vr_\infty > 0 \ \mbox{as}\ |x| \to \infty.
\]

\subsection{Constitutive relations}
\label{cr}

Keeping in mind our application to model the star dynamics, we impose several restrictions on the state equation.

\subsubsection{Equations of state}

Similarly to \cite{FeNo6A}, we suppose the EOS for the pressure takes the form
\begin{equation} \label{pressure0}
p(\vr, \vt) = p_m(\vr, \vt) + p_r(\vt),
\end{equation}
where
\begin{equation} \label{pressure}
p_m(\vr,\vt) =  \vt^{\frac{5}{2}} P(Z),\quad  p_r(\vt) = \frac{a}{3} \vt^4,\quad  a > 0,\quad  Z = \frac{\vr}{\vt^{\frac{3}{2}}}.
\end{equation}
The first component represents the general form of the EOS of monoatomic gas satisfying 
\[
p_m = \frac{2}{3} \vr e_m 
\]
compatible with Gibbs' equation \eqref{Gib}. The second component is the radiation pressure relevant in models in astrophysics, see e.g. Battaner 
\cite{BATT}.

Accordingly, the internal energy reads  
\begin{equation} \label{internal energy0} 
e(\vr, \vt) = e_m(\vr, \vt) + e_r (\vr, \vt),  
\end{equation}
where
\begin{equation} \label{internal energy} 
e_m(\vr, \vt) = 
\frac{3}{2} \frac{\vt^{\frac{5}{2}}}{\vr}P(Z),\quad   e_r(\vr, \vt) = \frac{a}{\vr} \vt^4. 
\end{equation}
In addition, in accordance with Gibbs' equation \eqref{Gib},
\begin{equation} \label{entropy0}
s(\vr, \vt) = s_m(\vr, \vt) + s_r (\vr, \vt),\  
	\end{equation}
 where 
 \begin{equation} \label{entropy}
s_m (\vr ,\vt) =  \mathcal{S}(Z),\quad   s_r(\vr, \vt) =  \frac{4a}{3} \frac{\vt^3}{\vr},	
	\end{equation}
and 
\begin{equation} \label{c5}
	\mathcal{S}'(Z) = -\frac{3}{2} \frac{ \frac{5}{3} P(Z) - P'(Z) Z }{Z^2}.
\end{equation}

Next, let us compute the specific heat at constant volume 
\[
c_v (\vr, \vt) = \frac{\partial e(\vr, \vt)}{\partial \vt} = \frac{9}{4} \frac{ \frac{5}{3} P(Z) - P'(Z) Z }{Z} + \frac{4a}{\vr} \vt^3.	 
\]
Adopting the thermodynamics stability hypothesis we get certain constraint to be satisfied by $P$, namely 
\begin{align}
\frac{\partial p(\vr, \vt)}{\partial \vr} &> 0 \ \Rightarrow \ P'(Z) > 0 \ \mbox{for}\ Z > 0, \br
c_v(\vr, \vt) &> 0 \ \Rightarrow \  \frac{ \frac{5}{3} P(Z) - P'(Z) Z }{Z} > 0 \ \mbox{for}\ Z > 0,
\label{c2}
\end{align}
cf. M\" uller and Ruggeri \cite[Chapter 4, Section 2.4]{MURU}.

We also impose the Third law of thermodynamics in the area of degenerate gas $Z >> 1$, 
\[
\lim_{\vt \to 0+} s(\vr, \vt) = 0 \ \mbox{for any}\ \vr > 0,
\]
or 
\begin{equation} \label{3rdlaw}
\mathcal{S}(Z) \to 0 \ \mbox{as}\ Z \to \infty,
\end{equation}
cf. Belgiorno \cite{BEL1}, \cite{BEL2}. In view  \eqref{3rdlaw}, we may strengthen \eqref{c2} to 
\begin{equation} \label{c2a}
0 < 	\frac{ \frac{5}{3} P(Z) - P'(Z) Z }{Z} < \Ov{P} \ \mbox{for}\ Z > 0. 
	\end{equation}  
 
Finally, we anticipate the EOS in the non--degenerate area of moderate values of $Z$ to obey the standard Boyle--Mariotte law corresponding 
to $P(Z) \approx Z$. This motivates the hypothesis
\begin{equation} \label{c2b}
P(Z) \in C^2[0, \infty),\quad P(0) = 0,\quad  P''(0) =  0 < \underline{P} \leq \frac{P(Z)}{Z} \ \mbox{for all}\ Z > 0.  
\end{equation}
Note that \eqref{c2b} implies 
\begin{equation} \label{c2ba}
	P'(0) > 0.
\end{equation}

\subsubsection{Transport coefficients}

Motivated by \cite{ChauFei}, \cite[Chapter 12]{FeiNovOpen}, we suppose
the transport coefficients $\mu$, $\eta$, and $\kappa$ in \eqref{NR} and \eqref{FL}, are continuously differentiable functions of the temperature $\vt$, satisfying:
\begin{align}
	0 < \underline{\mu} \left(1 + \vt \right) &\leq \mu(\vt) \leq \Ov{\mu} \left( 1 + \vt \right),\ 
	|\mu'(\vt)| \leq c \ \mbox{for all}\ \vt \geq 0, \label{mu}\\
	0 &\leq  \eta(\vt) \leq \Ov{\eta} \left( 1 + \vt \right),  \label{eta}\\
	0 < \underline{\kappa} \left(1 + \vt^\beta \right) &\leq  \kappa(\vt) \leq \Ov{\kappa} \left( 1 + \vt^\beta \right),\; \beta>6.
	\label{kappa}
\end{align}
The growth conditions imposed on the viscosity coefficients are standard and correspond to the so--called hard sphere gas model. The growth 
of the heat conductivity $\kappa$ is enforced by radiation, however, 
the restriction $\beta > 6$ is purely technical, cf. \cite{ChauFei}.

The paper is organized as follows. In Section \ref{w}, we introduce the concept of weak solution motivated by \cite{FeiNovOpen}. The main results concerning the 
existence of global--in--time weak solutions and the weak--strong uniqueness property are stated in Section \ref{M}. The proof of existence is given in Section \ref{P}, while the weak--strong uniqueness is established in Section \ref{WS}.

\section{Weak formulation}
\label{w}

Here and hereafter, we suppose $\Omega \subset \R^3$ is an exterior domain with 
regular boundary of class at least $C^{2 + \nu}$. Similarly, we suppose the boundary/far field temperature distribution is determined by a function $\vtB$,
\begin{equation} \label{bd2}
	\vtB \in C^{2 + \nu}([0,T] \times \R^3),\vtB \geq\vt_\infty \ \mbox{on}\ \R^3, 
	\ \vtB (t, x) = \vt_\infty \ \mbox{for}\ |x| > r,
\end{equation} 
where $r > 0$ is chosen so large that 
\[
\left\{ |x| > \frac{r}{2} \right\} \subset \Omega.
\]

Finally, in accordance with the no--slip boundary conditions \eqref{BC1}, we fix 
\begin{equation} \label{bd3}
	\vr = \Ov{\vr},\ \vu = 0 \ \mbox{in}\ \R^3 \setminus \Omega,
\end{equation}
where $\Ov{\vr} = \Ov{\vr}(x) \geq 0$ -- a bounded measurable function -- is the mass distribution of the ``hard core'' appearing in \eqref{G1}.

\subsection{Initial data}
 
A quantity that plays a crucial role in the analysis of 
stability of weak solutions is the	 
\emph{relative energy}
\[
E \left( \vr, \vt, \vu \ \Big| \tvr, \tvt, \tvu \right) = 
\frac{1}{2} \vr |\vu - \tvu |^2 + \vr e(\vr, \vt) - \tvt \vr s(\vr, \vt)- \Big( e(\tvr, \tvt) - \tvt s(\tvr, \tvt) + \frac{p(\tvr, \tvt)}{\tvr} \Big) \vr + p(\tvr, \tvt).
\]
As shown in \cite[Chapter 3, Section 3.1]{FeiNovOpen}, the relative energy 
represents a Bregman distance when expressed in terms of the conservative entropy variables $\vr, \vm = \vr \vu, S = \vr s(\vr, \vt)$.  Accordingly, we prescribe the initial state of the system in the form  
\begin{equation} \label{bd1}
	\vr(0, \cdot) = \vr_0,\ \vr_0 > 0,\ \vr \vu(0, \cdot) = \vm_0 = \vr_0 \vu_0,\ 
	\vr s(\vr, \vt)(0, \cdot) = S_0 = \vr_0 s(\vr_0, \vt_0).
\end{equation}

In addition, we suppose 
that the value of the relative energy evaluated with respect to the 
``far field data'' $\tvr = 0$, $\tvu = 0$, and
\[
\tvt \in C^1([0,T] \times \Ov{\Omega}),\ \tvt > 0, \ \tvt|_{\partial \Omega} =
 \vtB,\ \tvt \to \vt_\infty \ \mbox{as}\ |x| \to \infty 
\]
should be finite. As the expression 
\[
\Big( e(\tvr, \tvt) - \tvt s(\tvr, \tvt) + \frac{p(\tvr, \tvt)}{\tvr} \Big)
\]
is independent of the radiation component (cf. \eqref{pressure} -- \eqref{entropy}), a suitable hypothesis on the initial data reads
\[
\intO{ \left[ \frac{1}{2} \vr_0 |\vu_0|^2 + \vr_0 e(\vr_0, \vt_0) - \tvt 
\vr_0 s(\vr_0, \vt_0) + \frac{a}{3} \tvt^4 \right] } < \infty,	
\]
or, equivalently,  
\begin{equation} \label{bd1a}
\intO{ \left[ \frac{1}{2} \vr_0 |\vu_0|^2 + \vr_0 e(\vr_0, \vt_0) - \tvt 
	\vr_0 s(\vr_0, \vt_0) + \frac{a}{3} \vt^4_\infty \right]} < \infty.
\end{equation}
In accordance with hypothesis \eqref{entropy}, the expression 
$\vr_0 s(\vr_0, \vt_0)$ behaves like $\vr_0 \log(\vr_0)$ for $|x| \to \infty$, which imposes a decay 
rate on the initial distribution of the density
\begin{equation} \label{bd1b}
\vr_0 > 0,\ \intO{ \vr_0 }= m_0< \infty,\ \intO{ \vr_0 |\log \vr_0| } < \infty.	
	\end{equation} 

\subsection{Weak formulations of the field equations}

Similarly to \cite{ChauFei}, \cite[Chapter 12]{FeiNovOpen} the concept of weak solution to the system of equations \eqref{NSF1}--\eqref{NSF3} is based on replacing the internal energy balance \eqref{NSF3} by the entropy balance (inequality) \eqref{ent_eq} supplemented with the balance of ballistic energy 
evaluated at a reference temperature $\tvt$.

\subsubsection{Equation of continuity, mass conservation}

With our convention \eqref{bd3}, the equation of continuity  \eqref{NSF1} is 
satisfied on the whole space $\R^3$:
\begin{align}
	\int_0^T &\intR{ \Big[ \vr \partial_t \varphi + \vr \vu \cdot \Grad \varphi \Big]} \dt = - \intR{ \vr_0 \varphi(0, \cdot)  } 
	\label{w:ce}
\end{align}
for any $\varphi \in C^1_c([0,T) \times \R^3)$.

In addition, we impose the renormalized equation of continuity in the form
 \begin{align}
	\int_0^T &\intR{ \Big[ b(\vr) \partial_t \varphi + b(\vr) \vu \cdot \Grad \varphi -(b^\prime(\vr) \vr -b(\vr)\Div \vu \varphi )\Big]} \dt = - \intR{ b(\vr_0) \varphi (0, \cdot) } 
	\label{w:rce}
\end{align}
for any $\varphi \in C^1_c([0,T) \times \R^3)$, and $b \in C^1[0, \infty)$, $b' \in C_c[0, \infty)$.

\subsubsection{Momentum equation}

In view of the no-slip boundary condition \eqref{BC1} together with the far field condition \eqref{FFC2}, we require
\[
\vu \in L^2(0,T; D^{1,2}_0(\Omega;\R^3)), 
\]
where $D^{k,q}(\Omega)$ denotes space of locally integrable functions s.t. $\| \Grad^k \vc{v} \|_{L^q(\Omega;\R^3)}$ is bounded. In accordance with \eqref{bd3}   
we may assume 
\begin{equation} \label{m:class}
	\vu  \in L^2(0,T; D^{1,2}_0(\R^3;\R^3)) 
\end{equation} 
extending $\vu$ to be zero outside $\Omega$.  The weak formulation of the momentum 
equation \eqref{NSF2} reads
\begin{align}
	\int_0^T &\intR{ \Big[ \vr \vu \cdot \partial_t \bfphi + \vr \vu \otimes \vu : \Grad \bfphi + p(\vr, \vt) \Div \bfphi \Big] } \dt \br &= 
	\int_0^T \intR{ \Big[ \mathbb{S} (\vt, \Ds \vu) : \Ds \bfphi - \vr  \Grad G \cdot \bfphi \Big] } \dt - \intR{ \vr_0 \vu_0 \cdot \bfphi (0, \cdot) },   
	\label{w:me}
\end{align}
for any $\bfphi \in C^1_c([0,T) \times {{\R^3}};\R^3)$, where the gravitational 
potential $G$ solves the Poisson equation 
\begin{equation} \label{sg}
	- \Div (\mathfrak{g} \Grad G) = \vr \ \mbox{in}\ \R^3,\ 
	{\mbox{where we have set}}\ \vr = \Ov{\vr} \ \mbox{in}\ 
	\R^3 \setminus \Omega.
\end{equation}
\subsubsection{Entropy balance}

Similarly to \cite[Chapter 3]{FeNo6A}, we impose the entropy balance in the form
\begin{align}
	&- \intR{ \vr_0 s(\vr_0, \vt_0) \varphi (0, \cdot) } - \int_0^T \intR{ 
		\left[ \vr s \partial_t \varphi + \vr s \vu \cdot \Grad \varphi + \frac{\vc{q}}{\vt} \cdot \Grad \varphi \right] } \dt \br &\geq 
	\int_0^T \intR{ \frac{\varphi}{\vt} \left( \mathbb{S}(\vt, \Ds \vu ): \Ds \vu - \frac{\vc{q}(\vt, \Grad \vt) \cdot \Grad \vt }{\vt} \right) } \dt 
	\label{w:ei2}
\end{align}
for  any $\varphi \in C^1_c([0,T) \times \Omega)$, $\varphi \geq 0$. Here, in accordance with our convention 
\eqref{bd2}, we require 
\begin{equation} \label{bd4}
	\vt - \vtB \in L^2(0,T; W^{1,2}_0(\Omega)).
\end{equation}	

\begin{Remark} \label{RR1}
As $\vt$ satisfies strictly positive far field condition $\vt_\infty$, we 
impose $W^{1,2}_0(\Omega)$ in \eqref{bd4} instead of the ``expected'' weaker assumption $D_0^{1,2}(\Omega)$, cf. the ballistic energy inequality \eqref{w:bei} below.

	\end{Remark}

\subsubsection{Ballistic energy balance}
\label{wf:be}

Following the strategy of \cite{ChauFei}, \cite{FeiNovOpen}, we close the weak formulation of the problem by imposing a 
\emph{ballistic} energy balance. Consider a function 
\[
\tvt \in C^1([0,T] \times\Ov{\Omega}),\ \tvt = \vtB \ \mbox{on}\ \partial\Omega, \ 
\tvt > 0,\ \tvt - \vt_\infty \in C^1_c([0,T] \times\Ov{\Omega}).
\]
The ballistic energy inequality reads
\begin{align}  
	& \intO{ \left[ \frac{1}{2} \vr  |\vu  |^2 + \vr  e  - \tvt \vr  s + \frac{a}{3} \vt_\infty^4    \right] (\tau, \cdot) }\dt +
	\int_0^\tau  \intO{ \frac{\tvt }{\vt }	 \left( \mathbb{S}  : \Ds \vu  - \frac{\vc{q}  \cdot \Grad \vt  }{\vt } \right) }\dt  \br
	&\leq  \intO{ \left[ \frac{1}{2} \vr_0 |\vu_0  |^2 + \vr_0  e(\vr_0, \vt_0) - \tvt(0, \cdot)  \vr_0  s (\vr_0, \vt_0) + \frac{a}{3} \vt^4_\infty  \right] } + 	\int_0^\tau  \intO{ \vr \Grad G \cdot \vu  } \dt \br 
	&\ \ \ - 	\int_0^\tau  \intO{ \left[ \vr  s  \left( \partial_t \tvt  + \vu  \cdot \Grad \tvt  \right) + \frac{\vc{q} }{\vt} \cdot \Grad \tvt  \right] } \dt 
 \label{w:bei}
\end{align}
for a.a. $0 < \tau < T$, cf. \eqref{bd1a}. 

\subsection{Weak solutions}  

Having collected all preliminary material, we introduce the concept of 
weak solution to the NSF system. 

\begin{Definition}[{\bf Weak solution}] \label{DW1}
	
	Let the initial data belong to the class \eqref{bd1a}, \eqref{bd1b}, and let 
	the boundary/far field temperature be determined by \eqref{bd2}. 
We say that a trio $(\vr, \vt, \vu)$ is \emph{weak solution} of 
the NSF system \eqref{NSF1}--\eqref{FL}, with the boundary conditions \eqref{BC1}, \eqref{BC2}, and the total mass constraint \eqref{BC1a}, if the relations 
\eqref{w:ce}--\eqref{w:bei}, together with \eqref{BC1a}, are satisfied in $(0,T) \times \Omega$.
	
	\end{Definition}

\section{Main results}
\label{M}

We are ready to state our main results on the weak solutions to the NSF system.

\subsection{Global--in--time existence}

Our first result claims global--in--time existence of weak solution for arbitrary (large) initial/boundary data. 

\begin{mdframed}[style=MyFrame]

\begin{Theorem}[{\bf Global existence}] \label{Maintheorem}
	
	Let $\Omega \subset\R^3$ be an exterior domain of class $C^{2 + \nu}$. Suppose the pressure $p$, the internal energy $e$, and the entropy $s$ 
	comply with \eqref{pressure}--\eqref{c2b}, and the transport coefficients $\mu$, $\lambda$, and $\kappa$ satisfy \eqref{mu}--\eqref{kappa}.
	Let the boundary/far field data $\vtB$ be determined by
	a function
	\[
	\vtB \in C^{2 + \nu}([0,T]\times\R^3), \inf_{[0,T] \times \partial \Omega},   
	\inf_{[0,T] \times \partial \Omega} \vtB > \vt_\infty,\ 
	\vtB - \vt_\infty \in C^{2 + {\cblue{\nu}}} _c ([0,T] \times\R^3).
	\]
	Let the initial data belong to the class 
	\begin{align} \label{class}
	\vr_0 \in L^{\frac{5}{3}}(\Omega) \cap L^1(\Omega),\ \vr_0 \geq 0, \ \vr_0 \log(\vr_0) , \vr_0 |\vu_0|^2 \in L^1(\Omega),\br 
	\vt_0 \in L^\infty(\Omega),\ (\vt_0 - \vt_\infty) \in L^2(\Omega),\ \vt_0 > 0.
	\end{align}
	
	Then the NSF system \eqref{NSF1}--\eqref{FL}, coupled with the Poisson equation \eqref{G1}, and endowed with the boundary conditions \eqref{BC1}--\eqref{FFC2}
	admits a weak solution $(\vr, \vt, \vu)$ in $(0,T) \times \Omega$ in the sense specified in Definition \ref{DW1}.
	
	\end{Theorem}
\end{mdframed}
{Note that the initial data \eqref{class} satisfy the conditions \eqref{bd1a}, \eqref{bd1b}.}
Theorem \ref{Maintheorem} will be proved in Section \ref{P}.

\subsection{Weak--strong uniqueness}	

In order to formulate the weak--strong uniqueness result we have to specify the class of strong solutions. As already pointed out, this is a delicate issue as, to the best of our knowledge, the existence of strong solutions even local in time is not known for this particular choice of the boundary/far field conditions. We therefore anticipate the regularity available for similar problems. Following \cite{HuLiXi}, we consider strong solutions $(\tvr, \tvt, \tvu)$ belonging to the class: 
\begin{align}
	\tvr &\in C([0,T]; (W^{3,2} \cap L^1) (\Omega)),\ \tvr > 0 ,\  
	\tvr \log(\tvr) \in C([0,T]; L^1(\Omega)), \br
	(\tvt - \vtB) &\in C([0,T]; W^{1,2}_0 \cap D^{3,2}(\Omega)) \cap 
	L^2(0,T; D^{4,2}(\Omega)),\quad \inf_{(0,T) \times \Omega} \tvt > 0 \br
	\partial_t \tvt  &\in L^\infty(0,T; W^{1,2} (\Omega)) \cap 
	L^2(0,T; D^{2,2}(\Omega)), \br
	\tvu &\in C([0,T]; W^{1,2}_0 \cap D^{3,2}(\Omega;\R^3)) \cap 
	L^2(0,T; D^{4,2}(\Omega;\R^3)), \br
	\partial_t \tvu &\in L^\infty(0,T; D^{1,2}_0 (\Omega;\R^3)) \cap 
	L^2(0,T; D^{2,2}(\Omega;\R^3)), \br 
	\sqrt{\tvr}\partial^2_{t,t} \tvu &\in L^2((0,T) \times \Omega;\R^3).
		\label{classS}
	\end{align}

\begin{Remark}\label{PR1}
		Note that $\partial_t \tvt$ enjoys better integrability than $\partial_t \vu$ for $|x| \to \infty$ as the internal energy contains the radiative component 
		that is not multiplied by $\vr$.
		\end{Remark}

In addition to the above hypotheses, we need more explicit decay properties of 
the strong solution specified in Theorem \ref{TM2} below. They basically require the pressure to decay at least as fast as the density and are compatible with the polynomial decay of $\tvr$, 
\[
\tvr \approx |x|^{-\beta},\ \beta > 1 \ \mbox{as}\ |x| \to \infty.
\]		
For the sake of simplicity, we also asuume $G = 0$. 


	

	
	
	


\begin{mdframed}[style=MyFrame]

\begin{Theorem}[{\bf Weak--strong uniqueness}] \label{TM2}
	
Suppose the NSF system admits a 
strong solution $(\tvr, \tvt, \tvu)$ in $(0,T) \times \Omega$ belonging to the regularity class \eqref{classS}. Let the initial/boundary data be the same as in  Theorem \ref{Maintheorem}, and, in addition
\begin{equation} \label{jjj}
\Grad \log (\tvr_0 )\in L^\infty(\Omega, \R^3),\quad \tvr_0 \aleq \frac{1}{|x|^2}
\ \mbox{as}\ |x| \to \infty,
\end{equation}
Finally, suppose that 
$(\tvr, \tvt,\tvu)$ satisfies the decay property
 \begin{equation}\label{proper0}
 \frac{1}{\tvr} \Div \mathbb{S} (\tvt, \Grad \tvu) \aleq 1 \ \mbox{as}\ |x| \to \infty.
 \end{equation}
Then, $\tvr = \vr$, $\tvt = \vt$, and $\tvu = \vu$, where $(\vr, \vt, \vu)$ 
is a weak solution of the same problem in the sense of Definition \ref{DW1}. 	
	
	\end{Theorem}	

	\end{mdframed}

 Theorem \ref{TM2} is proved in Section \ref{WS}.

 \begin{Remark}
 Hypotheses \eqref{jjj}, \eqref{proper0} yield 
   	\begin{equation} \label{proper}
|\Grad \tvr| + |\Grad \tvt| \aleq \tvr \aleq \frac{1}{|x|^2}
\ \mbox{as}\ |x| \to \infty.
	\end{equation} 
Indeed the decay of the initial density $\vr_0$ persists at any positive time as long as the velocity is continuously differentiable (twice). This follows from the the equation of continuity. 
To see the decay of the temperature, we first rewrite \eqref{NSF2} for classical solutions, we obtain
\[
    4 \tvt^3 \frac{\nabla_x \tvt}{\tvr }= -\partial_t  \tvu - \tvu \cdot \nabla_x \tvu - \frac{1}{\tvr}\nabla_x p_m(\tvr,\tvt) + \frac{1}{\tvr} \Div \bb{S}(\tvt , \Ds \tvu) +  \vc{g}.
\]
The right hand side of the above is bounded thanks to hypothesis \eqref{proper0} and the observation above. In particular, seeing that $p_m(\tvr,\tvt)\approx\tvr\tvt + \tvr^\gamma $, we have 
\[ \frac{1}{\tvr}\nabla_x p_m (\tvr,\tvt)= \tvt \frac{\Grad \tvr}{\tvr} + \Grad \tvt+ \tvr^{\gamma-1}\frac{\Grad \tvr}{\tvr}, \]
which is bounded.

\end{Remark}

\section{Proof of Theorem \ref{Maintheorem}}
\label{P}	

The proof of Theorem \ref{Maintheorem} is based on uniform bounds for a sequence of suitable approximate solutions.

\subsection{Method of invading domains} 

We approximate our problem by considering it  on a family of bounded domains 
\begin{equation} \label{PP1}
	\Omega_n = \Omega \cap \left\{ x \in\R^3 \ \Big|\ |x| \leq r + n \right\},
\end{equation}
supplemented with the boundary conditions 
\begin{equation} \label{P3}
	\vun|_{|x| = {r}+ n} = 0,\ \vtn|_{|x| = r + n} = \vt_\infty. 
\end{equation}
The initial data for the approximate problems are given by restrictions of 
$(\vr_0, \vt_0, \vu_0)$ on $\Omega_n$. The weak formulation of the approximate problem is exactly the same as in Section \ref{w}, with $\Omega$ being replaced by $\Omega_n$.

The existence of weak solutions $(\vrn, \vtn, \vun)_{n=1}^\infty$ in $(0,T) \times \Omega_n$ was proved in \cite{ChauFei}, \cite{FeiNovOpen}, or, by a penalization method in \cite{BaFeLMMiYu}. Note that adding the Poisson equation \eqref{sg} does not bring any essential difficulty, see e.g. \cite{DFPS1}. Moreover, the limit 
$n \to \infty$ for problems with \emph{conservative} boundary/far field conditions was performed in several papers:  
\cite{DFPS1}, Jessl\' e, Jin, Novotn\' y, \cite{JeJiNo}, Poul \cite{Poul1}. In particular, the (non--trivial) issues concerning compactness of the 
sequence $(\vrn, \vtn, \vun)_{n=1}^\infty$ in the \emph{interior} of $\Omega$ are nowadays well understood. In view of these arguments, the proof 
of Theorem \ref{Maintheorem} reduces to showing suitable \emph{uniform} bounds on the sequence  $(\vrn, \vtn, \vun)_{n=1}^\infty$ and passing to the limit in the ballistic 
energy inequality \eqref{w:bei}.  

\subsection{Uniform bounds}

In accordance with the previous plan, our ultimate goal is to establish uniform 
bounds on the sequence of approximate solutions. 

\subsubsection{Mass conservation}

First observe that a total mass of the approximate solutions is a conserved quantity, 
\begin{equation} \label{P4}
	\intOn{ \vrn(t, \cdot) } = \intOn{ \vr_0 } = m_{0,n}\ \mbox{for any}\ t \in (0,T).
\end{equation}
As obviously 
\[
m_{0,n} \to m_0 = \intO{ \vr_0 },
\]
relation \eqref{P4} yields a uniform bound 
\begin{equation} \label{P5}
	\| \vrn(t, \cdot) \|_{L^1(\Omega_n)} \aleq 1 \ \mbox{independent of}\ n = 1,2, \dots
\end{equation}
Here and hereafter, the symbol $a \aleq b$ means there is a positive constant $C$ independent of $n$ such that $a \leq Cb$.

\subsubsection{Ballistic energy structural bounds}
We start from reportting the following result on the structural properties of EOS proved in \cite[Section 2.1.1]{FeiLuSun}: 

\begin{Lemma} \label{LP1}
Under the hypotheses of Section \ref{cr}, we have 
\begin{equation} \label{P1}
P(Z) = p_\infty Z^{\frac{5}{3}} + \widetilde{P}(Z),\ p_\infty > 0,	
	\end{equation}
and 
\begin{align} \label{P2}
\widetilde{P} \in C^1[0, \infty),\ \widetilde{P}(0) = 0,\ \widetilde{P}'(0) > 0, \ \widetilde{P} \geq 0, \ \frac{\widetilde{P}(Z)}{Z} \to 0 \ \mbox{as}\ Z \to \infty, \br |\widetilde{P}'(Z)| \leq c 
\ \mbox{for all}\ Z \geq 0. 
\end{align} 
\end{Lemma}

Now, we examine the ballistic energy associated to the far--field temperature, namely,  
\[
\frac{1}{2} \vr |\vu|^2 + \vr e(\vr, \vt) - \vt_\infty \vr s(\vr, \vt) = 
\frac{1}{2} \vr |\vu|^2 + \frac{3}{2} \frac{\vt^{\frac{5}{2}}}{\vr} P(Z) - \vt_\infty \vr \mathcal{S} \left( \frac{\vr}{\vt^{\frac{3}{2}}}  \right) + a \vt^4  - 
\frac{4a}{3} \vt_\infty \vt^3.
\]
It follows from the structural properties of the pressure EOS \eqref{c2}, \eqref{P1} that 
\begin{equation} \label{PP6a} 
\frac{3}{2} \frac{\vt^{\frac{5}{2}}}{\vr} P(Z) \geq c_1 \vr^{\frac{5}{3}} + c_2 \vr \vt ,\ c_1, c_2 > 0. 
\end{equation}
Moreover, by virtue of \eqref{entropy}--\eqref{3rdlaw}, 
\begin{equation} \label{PP6}
0 \leq \vr \mathcal{S} \left( \frac{\vr}{\vt^{\frac{3}{2}}}  \right) \aleq \Big( \vr + \vr |\log (\vr)| + \vr \log^+(\vt) \Big).
\end{equation}
Indeed, on the one hand, it follows from \eqref{c5}, \eqref{3rdlaw}, and \eqref{c2a} that 
\[
0 \leq \mathcal{S}\left( \frac{\vr}{\vt^{\frac{3}{2}}}  \right) \leq \mathcal{S}(1) \quad \mbox{if} \  \frac{\vr}{\vt^{\frac{3}{2}}} \geq 1.
\]
On the other hand, if $\frac{\vr}{\vt^{\frac{3}{2}}} < 1$, we have
\begin{align}
0 \leq \mathcal{S} \left( \frac{\vr}{\vt^{\frac{3}{2}}} \right) &= 
\mathcal{S}(1) - \int_{ \frac{\vr}{\vt^{\frac{3}{2}}} }^1 \mathcal{S}'(y) \ \D y
\leq \mathcal{S}(1) + \frac{3 \Ov{P}}{2} \int_{ \frac{\vr}{\vt^{\frac{3}{2}}} }^1 \frac{1}{y} \ \D y    = \mathcal{S}(1) - \frac{3 \Ov{P}}{2} \log \left( \frac{\vr}{\vt^{\frac{3}{2}}} \right) \br 
&= 
\mathcal{S}(1) - \frac{3 \Ov{P}}{2} \log (\vr) + \frac{9 \Ov{P}}{4} \log(\vt)
\leq \mathcal{S}(1) + \frac{9 \Ov{P}}{4} \Big( |\log(\vr)| + \log^+(\vt) \Big). 
\nonumber
\end{align}

Finally,
\[ 
a \vt^4  - 
\frac{4a}{3} \vt_\infty \vt^3 = H(\vt, \vt_\infty) - \frac{a}{3} \vt_\infty^4, 
\]
where 
\begin{align} 
 &H(\vt, \vt_\infty)\equiv a \vt^4  - \frac{4a}{3} \vt_\infty \vt^3 + \frac{a}{3} \vt_\infty^4 , \br 
&H(\vt, \vt_\infty)  \ageq |\vt - \vt_\infty|^2 \ \mbox{if}\ \vt_\infty / 2 < \vt < 2 \vt_\infty, \br
&H(\vt, \vt_\infty)  \ageq 1 \ \mbox{if}\ \vt \leq \vt_\infty / 2, \ H(\vt, \vt_\infty)  \ageq  \vt^4 
\ \mbox{if}\ \vt \geq 2 \vt_\infty.
\label{P6}
\end{align}

Summing up the previous observations, we conclude there exists a positive constant $K$ depending on $a$, $\vt_\infty$ and the structural properties of $P$ 
such that
\begin{equation} \label{P7}
\vr e(\vr, \vt) - \vt_\infty \vr s(\vr, \vt) + 
K \vr \Big( 1  + |\log(\vr)| \Big) + \frac{a}{3} \vt_\infty^4 \ageq
\Big( \vr^{\frac{5}{3}} + \vr \vt + \vr |\log(\vr)| + H(\vt, \vt_\infty) \Big).
\end{equation}
Note that the term proportional $\vr \log^+(\vt)$ can be controlled 
via \eqref{PP6a} and $K \vr (1 + |\log(\vr)|)$ provided $K$ is large enough. 

\subsubsection{Controlling small values of $\vr$}

Our next goal is to control the logarithmic term $\vr|\log(\vr)|$ for $\vr$ approaching zero. To this end,  
consider the cut--off function 
\[
T_k(\vr) = \min \{ \vr, k \}.
\]
The renormalized equation \eqref{w:rce} yields
\begin{equation}\label{AA1}
\partial_t (\vr L_k(\vr)) + 
\Div (\vr L_k(\vr) \vu) + T_k (\vr) \Div \vu = 0,
\end{equation}
where
\[
L_k(\vr) = \int_1^\vr \frac{T_k(z)}{z^2} \ {\rm{d}} z.
\]
Observe that $\vr L_k(\vr) = \vr \log (\vr)$ for $\vr \leq k$. Note that $\vr \mapsto \vr L_k (\vr)$ is an affine function for 
$\vr \geq k$ so \eqref{AA1} can be justified by means of \eqref{w:ce}, \eqref{w:rce}.  

Integrating the renormalized equation of continuity, we get 
\begin{equation} \label{PP9}
\frac{\D }{\dt} \intOn{ \vr L_k(\vr) } = - 
\intOn{ T_k(\vr) \Div \vu},
\end{equation}
where 
\begin{align}
\left| \intOn{ T_k(\vr) \Div \vu} \right| \leq 
\| T_k(\vr) \|_{L^2(\Omega_n)} \| \Div \vu \|_{L^2(\Omega_n)}, 
\br
\| T_k(\vr) \|_{L^2(\Omega_n)}^2 \leq k \| \vr \|_{L^1(\Omega_n)}  \leq k m_0.
\label{PP9a}
\end{align}

\subsubsection{Ballistic energy inequality}

Finally, we examine the ballistic energy inequality \eqref{w:bei} reading now 
\begin{align}   
	& \intOn{ \left[ \frac{1}{2} \vrn  |\vun  |^2 + \vrn  e(\vrn, \vtn)  - \tvt \vr_n   s(\vr_n, \vt_n) + \frac{a}{3} \vt^4_\infty  \right] (\tau, \cdot) }  \br &+
	\int_0^\tau  \intOn{ \frac{\tvt }{\vt_n }	 \left( \mathbb{S}(\vt_n, \Ds \vu_n)  : \Ds \vu_n  + \frac{\kappa(\vt_n) | \Grad \vtn |^2  }{\vtn } \right) }\dt  \br
	&\leq  \intOn{ \left[ \frac{1}{2} \vr_0 |\vu_0  |^2 + \vr_0  e(\vr_0, \vt_0) - \tvt(0, \cdot)  \vr_0  s (\vr_0, \vt_0) + \frac{a}{3} \vt^4_\infty  \right] } + 	\int_0^\tau  \intOn{ \vrn \Grad G_n \cdot \vun  } \dt \br 
	&\ \ \ - 	\int_0^\tau  \intOn{ \left[ \vrn  s(\vrn, \vtn)  \left( \partial_t \tvt  + \vun  \cdot \Grad \tvt  \right) - \frac{\kappa (\vtn) \Grad \vtn }{\vtn} \cdot \Grad \tvt  \right] } \dt  \label{w:bein1}.
\end{align}
The first observation is that the class of functions $\tvt$ may be extended to 
\begin{equation} \label{P9}
\tvt \in W^{1, \infty}([0,T] \times \Omega_n),\ \tvt|_{\partial \Omega} = \vtB,\ \tvt|_{|x| = r+ n} = \vt_\infty.	
	\end{equation}
Moreover, it is easy to see that 
\begin{align}   
	& \intOn{ \left[ \frac{1}{2} \vrn  |\vun  |^2 + \vrn  e(\vrn, \vtn)  - \tvt \vr_n   s(\vr_n, \vt_n) + \frac{a}{3} \tvt^4  \right] (\tau, \cdot) }  \br &+
	\int_0^\tau  \intOn{ \frac{\tvt }{\vt_n }	 \left( \mathbb{S}(\vt_n, \Ds \vu_n)  : \Ds \vu_n  + \frac{\kappa(\vt_n) | \Grad \vtn |^2  }{\vtn } \right) }\dt  \br
	&\leq  \intOn{ \left[ \frac{1}{2} \vr_0 |\vu_0  |^2 + \vr_0  e(\vr_0, \vt_0) - \tvt(0, \cdot)  \vr_0  s (\vr_0, \vt_0) + \frac{a}{3} \tvt^4(0, \cdot)  \right] } + 	\int_0^\tau  \intOn{ \vrn g \Grad G_n \cdot \vun  } \dt \br 
	&\ \ \ - 	\int_0^\tau  \intOn{ \left[ \vrn  s(\vrn, \vtn)  \left( \partial_t \tvt  + \vun  \cdot \Grad \tvt  \right) - \frac{\kappa (\vtn) \Grad \vtn }{\vtn} \cdot \Grad \tvt  \right] } \dt + \int_0^\tau \intOn{ \frac{4a}{3} \tvt^3 \partial_t \tvt } \dt.  \label{w:bein}
\end{align}	
	
Keeping in mind 
\[
\left\{ |x| > \frac{r}{2} \right\} \subset \Omega,
\]	
we consider $\tvt = \tvt_r$ in the form 
\begin{equation} \label{P10}
\Del \tvt_r = 0 \ \mbox{in}\ \Omega \cap \{ |x| < r \},\ \tvt|_{\partial \Omega} = \vtB,\ \tvt_r |_{|x| = r} = \vt_\infty,\ 
\tvt_r = \vt_\infty \ \mbox{for}\ |x| > r.
\end{equation}	

Consequently, the inequality \eqref{w:bein} reduces to 
\begin{align}   
	& \int_{\Omega \cap \{ |x| < r \} } \left[ \frac{1}{2} \vrn  |\vun  |^2 + \vrn  e(\vrn, \vtn)  - \tvt_r \vr_n   s(\vr_n, \vt_n) + \frac{a}{3} \tvt_r^4  \right] (\tau, \cdot) \dx  \br
	&+ \int_{\Omega_n \cap \{ |x| \geq r \} } \left[ \frac{1}{2} \vrn  |\vun  |^2 + \vrn  e(\vrn, \vtn)  - \vt_\infty \vr_n   s(\vr_n, \vt_n) + \frac{a}{3} \vt_\infty^4  \right] (\tau, \cdot) \dx  \br
	 &+
	\int_0^\tau  \intOn{ \frac{\tvt_r }{\vt_n }	 \left( \mathbb{S}(\vt_n, \Ds \vu_n)  : \Ds \vu_n  + \frac{\kappa(\vt_n) | \Grad \vtn |^2  }{\vtn } \right) }\dt  \br
	&\leq  \intOn{ \left[ \frac{1}{2} \vr_0 |\vu_0  |^2 + \vr_0  e(\vr_0, \vt_0) - \tvt_r(0, \cdot)  \vr_0  s (\vr_0, \vt_0) + \frac{a}{3} \tvt_r^4 (0, \cdot) \right] } + 	\int_0^\tau  \intOn{ \vrn  \Grad G_n \cdot \vun  } \dt \br 
	&\quad - 	\int_0^\tau  \int_{\Omega \cap \{ |x| < r \} } \left[ \vrn  s(\vrn, \vtn)  \left( \partial_t \tvt_r  + \vun  \cdot \Grad \tvt_r  \right) - \frac{\kappa (\vtn) \Grad \vtn }{\vtn} \cdot \Grad \tvt_r  \right] \ \dx \dt \br 
	&\quad + \int_0^\tau \int_{\Omega \cap \{ |x| < r \} } \frac{4a}{3} \tvt^3_r \partial_t \tvt_r \dx \dt.   \label{P11}
\end{align}
The advantage of \eqref{P11} is that all integrals including $\tvt_r$ on the right--hand side are over a fixed bounded domain $\Omega \cap \{ |x| < r \}$, and consequently, they can be controlled in the same way as in the bounded domain case examined in \cite{ChauFei}, \cite{FeiNovOpen}.

Obviously, by virtue of \eqref{P10}, 
\[
\left| \int_0^\tau \int_{\Omega \cap \{ |x| < r \} } \frac{4a}{3} \tvt^3_r \partial_t \tvt_r \dx \dt \right| \leq C(r, T, \vtB).
\]

Next, in view of \eqref{P10}, 
\begin{align} 
\int_0^\tau &\int_{\Omega \cap \{ |x| < r \} } \frac{\kappa (\vtn) \Grad \vtn }{\vtn} \cdot \Grad \tvt_r   \ \dx \dt \br &= 
\int_0^\tau \int_{\partial \Omega } \mathcal{K} (\vtB) \Grad \tvt_r \cdot \vc{n} \ \D \sigma_x \dt + 
\int_0^\tau \int_{|x| = r } \mathcal{K} (\vtn) \Grad \tvt_r \cdot \vc{n} \ \D \sigma_x \dt,
\label{P12}
\end{align}
with $\mathcal{K}'(z) = \frac{\kappa(z)}{z}$, 
where 
\[
\left| \int_0^\tau \int_{\partial \Omega } \mathcal{K} (\vtB) \Grad \tvt_r \cdot \vc{n} \ \D \sigma_x \dt  \right| \leq C(r,T, \vtB).
\]
Moreover, as $\vtB > \vt_\infty$, the maximum principle yields $\tvt_r \geq \vt_\infty$ in $\Omega \cap \{ |x| < r \}$, in particular, 
\[
\Grad \tvt_r \cdot \vc{n}|_{|x| ={\cblue{r}}} \leq 0 \ \Rightarrow \  
\int_0^\tau \int_{|x| = r } \mathcal{K} (\vtn) \Grad \tvt_r \cdot \vc{n} \ \D \sigma_x \dt \leq 0.
\]

Summing up the previous observations we may infer that 	
\begin{align}   
	& \int_{\Omega \cap \{ |x| < r \} } \left[ \frac{1}{2} \vrn  |\vun  |^2 + \vrn  e(\vrn, \vtn)  - \tvt_r \vr_n   s(\vr_n, \vt_n) + \frac{a}{3} \tvt^4_r   \right] (\tau, \cdot) \dx  \br
	&+ \int_{\Omega_n \cap \{ |x| \geq r \} } \left[ \frac{1}{2} \vrn  |\vun  |^2 + \vrn  e(\vrn, \vtn)  - \vt_\infty \vr_n   s(\vr_n, \vt_n) + \frac{a}{3} \vt_\infty^4  \right] (\tau, \cdot) \dx  \br
	&+
	\int_0^\tau  \intOn{ \frac{\vt_\infty }{\vt_n }	 \left( \mathbb{S}(\vt_n, \Ds \vu_n)  : \Ds \vu_n  + \frac{\kappa(\vt_n) | \Grad \vtn |^2  }{\vtn } \right) }\dt  \br
	&\leq  \int_{\Omega \cap \{ |x| < r \} } \left[ \frac{1}{2} \vr_0 |\vu_0  |^2 + \vr_0  e(\vr_0, \vt_0) - \tvt_r(0, \cdot)  \vr_0  s (\vr_0, \vt_0) + \frac{a}{3} \tvt^4_r(0, \cdot)  \right] \dx \br  
	&\quad + \int_{\Omega \cap \{ |x| \geq r \} } \left[ \frac{1}{2} \vr_0 |\vu_0  |^2 + \vr_0  e(\vr_0, \vt_0) - \vt_\infty  \vr_0  s (\vr_0, \vt_0) + \frac{a}{3} \vt_\infty^4  \right] \dx
	\br &\quad + 	\int_0^\tau  \intOn{ \vrn \Grad G_n \cdot \vun  } \dt \br 
	&\quad - 	\int_0^\tau  \int_{\Omega \cap \{ |x| < r \} } \left[ \vrn  s(\vrn, \vtn)  \left( \partial_t \tvt_r  + \vun  \cdot \Grad \tvt_r  \right)   \right] \ \dx \dt  + C(r, T, \vtB).    \label{P13}
\end{align}

Next, by means of Korn's inequality and the hypotheses \eqref{mu}--\eqref{kappa} imposed on the transport coefficients, we deduce
\begin{align}
\int_0^\tau  &\intOn{ \left(  |\Grad \vu_n|^2  +  |\Grad \log(\vtn) |^2 + | \Grad \vtn^{\frac{\beta}{2}} |^2   \right) }\dt \br &+ \int_0^\tau \int_{\Omega \cap \{ |x| < r \} } \Big( | \vu_n |^2 + |\log(\vt_n)|^2 + |\vtn^{\frac{\beta}{2}} |^2   \Big) \dx \br
&\aleq 	\int_0^\tau  \intOn{ \frac{\vt_\infty }{\vt_n }	 \left( \mathbb{S}(\vt_n, \Ds \vu_n)  : \Ds \vu_n  + \frac{\kappa(\vt_n) | \Grad \vtn |^2  }{\vtn } \right) }\dt +  C(r, T, \vtB). 
\label{P13a}   
\end{align}
Thus going back to \eqref{P13} we get
\begin{align}   
	& \int_{\Omega \cap \{ |x| < r \} } \left[ \frac{1}{2} \vrn  |\vun  |^2 + \vrn  e(\vrn, \vtn)  - \tvt_r \vr_n   s(\vr_n, \vt_n) + \frac{a}{3} \tvt^4_r   \right] (\tau, \cdot) \dx  \br
	&+ \int_{\Omega_n \cap \{ |x| \geq r \} } \left[ \frac{1}{2} \vrn  |\vun  |^2 + \vrn  e(\vrn, \vtn)  - \vt_\infty \vr_n   s(\vr_n, \vt_n) + \frac{a}{3} \vt_\infty^4  \right] (\tau, \cdot) \dx  \br
	&+\int_0^\tau  \intOn{ \left(  |\Grad \vu_n|^2  +  |\Grad \log(\vtn) |^2 + | \Grad \vtn^{\frac{\beta}{2}} |^2   \right) }\dt \br &+ \int_0^\tau \int_{\Omega \cap \{ |x| < r \} } \Big( | \vu_n |^2 + |\log(\vt_n)|^2 + |\vtn^{\frac{\beta}{2}} |^2   \Big) \dx \br
	&\aleq  \int_{\Omega \cap \{ |x| < r \} } \left[ \frac{1}{2} \vr_0 |\vu_0  |^2 + \vr_0  e(\vr_0, \vt_0) - \tvt_r(0, \cdot)  \vr_0  s (\vr_0, \vt_0) + \frac{a}{3} \tvt^4_r(0, \cdot)  \right] \dx \br  
	&\quad + \int_{\Omega \cap \{ |x| \geq r \} } \left[ \frac{1}{2} \vr_0 |\vu_0  |^2 + \vr_0  e(\vr_0, \vt_0) - \vt_\infty  \vr_0  s (\vr_0, \vt_0) + \frac{a}{3} \vt_\infty^4  \right] \dx
	\br &\quad + \left|	\int_0^\tau  \intOn{ \vrn \Grad G_n \cdot \vun  } \dt \right| \br 
	&\quad + \left| 	\int_0^\tau  \int_{\Omega \cap \{ |x| < r \} } \left[ \vrn  s(\vrn, \vtn)  \left( \partial_t \tvt_r  + \vun  \cdot \Grad \tvt_r  \right)   \right] \ \dx \dt \right|  + C(r, T, \vtB).    \label{P14}
\end{align}
 
To exploit inequality \eqref{P7}, we combine the total mass conservation 
\eqref{P4}, \eqref{P5}, with the bounds \eqref{PP9}, \eqref{PP9a} obtaining 
\begin{equation} \label{P15}
\left[ \intOn{ \beta \vr - \vr L_1 (\vr) } \right]_{t = 0}^{t = \tau} \leq 
\delta \int_0^\tau \intOn{ |\Div \vu_n |^2 } + c(\delta, T, m_0) 
\end{equation}
for any $\beta > 0$, $\delta > 0$.We fix $\beta > 0$ large enough so that 
\[
\beta \vr - \vr L_1(\vr) \geq \vr \left( \frac{\beta}{2} + |\log^-(\vr)| \right).
\]
Now, we add \eqref{P15} to \eqref{P14} and use the coercivity estimates 
\eqref{P7} to conclude 
\begin{align}   
	& \intOn{ \left[ \vrn  |\vun  |^2 + \vrn \vtn + \vrn^{\frac{5}{3}} + 
		\vrn |\log(\vr_n)| + H(\vt_n, \vt_\infty)   \right] (\tau, \cdot) }  \br
	&+ \int_{\Omega_n \cap \{ |x| \geq r \} } \left[ \frac{1}{2} \vrn  |\vun  |^2 + \vrn  e(\vrn, \vtn)  - \vt_\infty \vr_n   s(\vr_n, \vt_n) + \frac{a}{3} \vt_\infty^4  \right] (\tau, \cdot) \dx  \br
	&+\int_0^\tau  \intOn{ \left(  |\Grad \vu_n|^2  +  |\Grad \log(\vtn) |^2 + | \Grad \vtn^{\frac{\beta}{2}} |^2   \right) }\dt \br &+ \int_0^\tau \int_{\Omega \cap \{ |x| < r \} } \Big( | \vu_n |^2 + |\log(\vt_n)|^2 + |\vtn^{\frac{\beta}{2}} |^2   \Big) \dx \br
	&\aleq  \intO{ \left[ \vr_0 |\vu_0  |^2 + \vr_0^{\frac{5}{3}} 
		+ \vr_0 |\log(\vr_0)| + |\vt_0 - \vt_\infty|^2 \right] } + \| \vt_0 \|_{L^\infty(\Omega)}^4 
		\br &\quad + \left|	\int_0^\tau  \intOn{ \vrn \Grad G_n \cdot \vun  } \dt \right| \br 
	&\quad + \left| 	\int_0^\tau  \int_{\Omega \cap \{ |x| < r \} } \left[ \vrn  s(\vrn, \vtn)  \left( \partial_t \tvt_r  + \vun  \cdot \Grad \tvt_r  \right)   \right] \ \dx \dt \right|  + C(r, T, \vtB).    \label{P16}
\end{align}

The integral on the right--hand side of \eqref{P16} containing the gravitational potential can be estimated in the same way as in \cite{DFPS1}. 

Finally, as the last integral in \eqref{P16} is evaluated over a fixed \emph{bounded set} $\Omega \cap \{ |x| < r \}$ it can be absorbed by the left--hand side of \eqref{P16} exactly as in \cite[Section 4.1]{ChauFei}
or \cite[Chapter 12, Section 12.4]{FeiNovOpen}. 


\subsubsection{Conclusion}

Applying Gronwall's lemma to \eqref{P16} and using regularity of the initial data stated in \eqref{class} we obtain the following bounds independent of $n=1,2,\dots$:
\begin{align} 
	\| \vr_n \|_{L^\infty(0,T; L^1 \cap L^{\frac{5}{3}}(\Omega_n))} &\aleq 1, \br
	\| \vrn \log(\vrn) \| _{L^\infty(0,T; L^1 (\Omega_n))} &\aleq 1 \br
	\| \sqrt{\vr_n} \vun \|_{L^\infty(0,T; L^1 \cap L^{2}(\Omega_n;R^3))} &\aleq 1, \br 
	\| H(\vtn, \vt_\infty) \|_{L^\infty(0,T; {{L^{1}}}(\Omega_n))} &\aleq 1, \br 
	\int_0^T \intOn{ \left(  |\Grad \vu_n|^2  +  |\Grad \log(\vtn) |^2 + | \Grad \vtn^{\frac{\beta}{2}} |^2   \right) }\dt &\aleq 1, \br 
	\int_0^T \int_{\Omega \cap K } \Big( | \vu_n |^2 + |\log(\vt_n)|^2 + |\vtn^{\frac{\beta}{2}} |^2   \Big) \ \dx \dt &\aleq 1
	\ \mbox{for any compact}\ K \subset\R^3.
	\label{P17}
	\end{align}

With the uniform bounds \eqref{P17} at hand and given the local character of the 
field equations, we can pass to the limit for $n \to \infty$ exactly as in the bounded domain case treated in detail in \cite{ChauFei}, and \cite[Chapter 12]{FeiNovOpen}. Thus we have completed the proof of Theorem \ref{Maintheorem}. $\Box$ 

\section{Weak-Strong uniqueness -- proof of Theorem \ref{TM2}}
\label{WS}

We proceed in several steps starting with an extension of the ballistic energy inequality. Recall that, in this part, we only consider the case $G=0$.

\subsection{Ballistic energy inequality}

First we show that validity of the ballistic energy inequality \eqref{w:bei}
can be extended to functions $\tvt$ belonging to the class \eqref{classS}. More specifically, we justify the inequality
\begin{align}  
	& \intO{ \left[ \frac{1}{2} \vr  |\vu  |^2 + \vr  e  - \tvt \vr  s + \frac{a}{3} \tvt^4    \right] (\tau, \cdot) }\dt +
	\int_0^\tau  \intO{ \frac{\tvt }{\vt }	 \left( \mathbb{S}  : \Ds \vu  - \frac{\vc{q}  \cdot \Grad \vt  }{\vt } \right) }\dt  \br
	&\leq  \intO{ \left[ \frac{1}{2} \vr_0 |\vu_0  |^2 + \vr_0  e(\vr_0, \vt_0) - \tvt(0, \cdot)  \vr_0  s (\vr_0, \vt_0) + \frac{a}{3} \tvt^4(0, \cdot)  \right] }  \br 
	&\ \ \ - 	\int_0^\tau  \intO{ \left[ \vr  s  \left( \partial_t \tvt  + \vu  \cdot \Grad \tvt  \right) + \frac{\vc{q} }{\vt} \cdot \Grad \tvt  \right] } \dt + \int_0^\tau \intO{ \frac{a}{3} \partial_t \tvt^4  } \dt .
	\label{WS1}
\end{align}
for a.a. $0 < \tau < T$. Note that \eqref{WS1} holds for functions $\tvt$ in the class specified in Section \ref{wf:be}, in particular $\partial_t \tvt$ is compactly supported and the extra integral on the right--hand side of \eqref{WS1} is well defined. Thus, if $\tvt$ is as in \eqref{classS}
we have

Next, seeing that 
\[
\intO{ \left(\frac{a}{3} \partial_t \tvt^4 - \vr s_r (\vr, \vt) \partial_t \tvt \right) } = \frac{4a}{3} \intO{ \left( \tvt^3 - \vt^3 \right) \partial_t \tvt } ,
\] 
where, by virtue of \eqref{classS} and \eqref{P17}, both 
$\left( \tvt^3 - \vt^3 \right)$ and $\partial_t \tvt$ belong to the space 
$L^2(0,T; L^2(\Omega))$. Consequently, we can extend \eqref{WS1} to all $\tvt$ satisfying 
\eqref{classS} by a density argument. In particular, \eqref{WS1} is satisfied as soon as $\tvt$ is the temperature component of a strong solution of the NSF system. 

\subsection{Relative energy inequality}

As already pointed out, the relative energy inequality associated to the 
NSF system reads 
\begin{equation}
	E \left(\vr, \vt, \vu \Big| \tvr, \tvt, \tvu \right) = 
	\frac{1}{2}\vr |\vu - \tvu|^2 + \vr e - \tvt \vr s - \Big( e(\tvr, \tvt) - \tvt s(\tvr, \tvt) + \frac{p(\tvr, \tvt)}{\tvr} \Big) \vr + p(\tvr, \tvt).
	\label{ws1}
\end{equation}
The strategy of the proof of Theorem \ref{TM2} is to describe the time evolution 
of $E \left(\vr, \vt, \vu \Big| \tvr, \tvt, \tvu \right)$, where 
$(\tvr, \tvt, \tvu)$ is the strong solution of the NSF system associated to the same data. Our starting point is the ballistic energy inequality \eqref{WS1}.

First, we observe that $\tvu$ can be used as a test function in the momentum equation \eqref{w:me} to obtain
\begin{align}
	\int_0^\tau &\intO{ \Big[ \vr \vu \cdot \partial_t \tvu + \vr \vu \otimes \vu : \Grad \tvu + p(\vr, \vt) \Div \tvu \Big] } \dt \br &= 
	\int_0^\tau \intO{ \mathbb{S} (\vt, \Ds \vu) : \Ds \tvu } \dt + \left[ \intR{ \vr \vu \cdot \tvu  } \right]_{t = 0}^{t=\tau}.   
	\label{WS2}
\end{align} 
Although $\tvu$ is not compactly supported, the above equality can be justified by a density argument.  For example, for any compactly supported $\bfphi$ we have 
\begin{align*}
-&\int_0^\tau \intO{ \Grad\vt^4 : \bfphi } \dt=\int_0^\tau \intO{ \vt^4 \Div \bfphi } \dt \\
&= 
\int_0^\tau \intO{ \left(\vt^4 - 4 \vt_\infty^3 (\vt - \vt_\infty) - \vt^4_\infty  \right)\Div \bfphi } \dt + \int_0^\tau \intO{ 4 \vt_\infty^3 (\vt - \vt_\infty) \Div \bfphi } \dt,
\nonumber 
\end{align*}
where 
\[
\left(\vt^4 - 4 \vt_\infty^3 (\vt - \vt_\infty) - \vt^4_\infty  \right) \in 
L^\infty (0,T; L^1(\Omega)), \quad (\vt - \vt_\infty) \in L^\infty (0,T; L^2(\Omega)). 
\]
Hence using a density argument, we obtain that the integral $\int_0^\tau \intO{ \vt^4 \Div \tilde\vu }$ is well defined provided that 
$\Div \tilde \vu \in L^\infty(0,T; BC \cap L^2(\Omega))$ in accordance with the regularity assumptions \eqref{classS}.

Similarly, we use $\frac{1}{2}|\tvu|^2$ as a test function in the equation of continuity \eqref{w:ce} obtaining 
\begin{align}
\left[ \intO{ \frac{1}{2} \vr |\tvu|^2 } \right]_{t=0}^{t=\tau}	= 
	\int_0^\tau \intO{ \Big[ \vr \tvu \cdot \partial_t \tvu + \vr (\vu \otimes \tvu) : \Grad \tvu  \Big]} \dt . 
	\label{WS3}
\end{align}

Adding \eqref{WS2} and \eqref{WS3} to \eqref{WS1} we get 
\begin{align}  
	& \intO{ \left[ \frac{1}{2} \vr  |\vu - \tvu  |^2 + \vr  e  - \tvt \vr  s + p_r (\tvt)    \right] (\tau, \cdot) }\dt +
	\int_0^\tau  \intO{ \frac{\tvt }{\vt }	 \left( \mathbb{S}  : \Ds \vu  - \frac{\vc{q}  \cdot \Grad \vt  }{\vt } \right) }\dt  \br
	&\leq  \intO{ \left[ \frac{1}{2} \vr_0 |\vu_0 - \tvu(0, \cdot)  |^2 + \vr_0  e(\vr_0, \vt_0) - \tvt(0, \cdot)  \vr_0  s (\vr_0, \vt_0) + p_r(\tvt) (0, \cdot)  \right] }  \br 
	&\quad  - 	\int_0^\tau  \intO{ \left[ \vr  s  \left( \partial_t \tvt  + \vu  \cdot \Grad \tvt  \right) + \frac{\vc{q} }{\vt} \cdot \Grad \tvt  \right] } \dt + \int_0^\tau \intO{ \partial_t p_r(\tvt)  } \dt \br 
	&\quad + \int_0^\tau \intO{ \Big[ \vr (\tvu - \vu) \cdot \partial_t \tvu + \vr \vu \otimes (\tvu - \vu): \cdot \Grad \tvu  \Big]} \dt \br 
	&\quad - \int_0^\tau \intO{ p(\vr, \vt) \Div \tvu  } \dt + \int_0^\tau \intO{ \mathbb{S} (\vt, \Ds \vu) : \Ds \tvu } \dt.
	\label{WS4}
\end{align}

\subsection{Singular terms}

To complete the proof of relative energy inequality, we have to add integrals describing the time evolution of
\begin{align}
\Big( e(\tvr, \tvt) - \tvt s(\tvr, \tvt) + \frac{p(\tvr, \tvt)}{\tvr} \Big) \vr + p_m(\tvr, \tvt) &= \Big( e_m(\tvr, \tvt) - \tvt s_m(\tvr, \tvt) + \frac{p_m(\tvr, \tvt)}{\tvr} \Big) \vr + p_m(\tvr, \tvt) \br
&= \Big( {\frac{5}{2}} e_m(\tvr, \tvt) - \tvt s_m(\tvr, \tvt) \Big) \vr + p_m(\tvr, \tvt).
\nonumber 
\end{align}
It turns out that the only problematic term is 
\[
\tvt s_m (\tvr, \tvt) \vr = \tvt \mathcal{S} \left( \frac{\tvr}{\tvt^{\frac{3}{2}}} 
	\right) \vr \approx \log(\tvr) \vr \ \mbox{for}\ |x| \to \infty, 
\]
where $\log(\tvr) \to {{-\infty}} \ \mbox{as}\ |x| \to \infty$.

Our goal is to show integrability of $- \vr \log(\tvr)$. Using the weak formulation of the equation of continuity \eqref{w:ce} we get, formally, 
\begin{equation} \label{WS5}
\left[ \intO{ - \vr \log(\tvr) } \right]_{t=0}^{t = \tau} =
- \int_0^\tau \intO{ \Big[ \vr \partial_t \log(\tvr) + \vr \vu \cdot \Grad \log(\tvr) \Big] } \dt  
\end{equation}
If $\tvr$ is a strong solution emanating from the same initial data $\vr_0$, we have 
\[
\left| \intO{ \vr_0 \log(\tvr)(0, \cdot) } \right| = \left| \intO{ \vr_0 \log(\vr_0) } \right| < \infty. 
\]
Thus we get $\vr \log(\tvr) \in L^\infty(0,T; L^1(\Omega))$ as soon as we control 
the integral on the right--hand side of \eqref{WS5}. As $\tvr > 0$ we have 
\[
\partial_t \log(\tvr) + \tvu \cdot \Grad \log(\tvr) = - \Div \tvu
\] 
and 
\[
\partial_t \left[ \partial_{x_j} \log(\tvr) \right] + 
\tvu \cdot \Grad \left[ \partial_{x_j} \log(\tvr) \right] = 
\partial_{x_j} \Div \tvu - \partial_{x_j} \tvu \cdot \Grad \log(\tvr),\ j=1,2,3. 
\]
As 
\begin{align}
\Grad \tvu &\in BC ([0,T] \times \Ov{\Omega}; \R^{3 \times 3}), \br
\tvu &\in C([0,T]; D^{3,2}(\Omega;\R^3)) \cap L^2(0,T; D^{4,2}(\Omega;\R^3))
\nonumber 
\end{align}
we deduce 
\begin{equation} \label{WS6}
\Grad \log(\tvr) = \frac{\Grad \tvr}{\tvr} \in L^\infty((0,T) \times \Omega)
\end{equation}
as long as $\frac{\Grad \tvr_0}{\tvr_0} \in L^\infty(\Omega)$ as required in 
Theorem \ref{TM2}. Seeing that 
\[
\partial_t \log(\tvr) = - \Div\tvu - \tvu \Grad \log(\tvr)
\]
we go back to \eqref{WS5} to conclude 
\begin{equation} \label{WS7}
{\rm ess} \sup_{t \in (0,T)} \| \vr \log (\tvr) \|_{L^1(\Omega)} \aleq 1.
\end{equation}

The estimates \eqref{WS6}, \eqref{WS7} allow us to evaluate the time 
evolution of the integral
\[ 
\intO{ \Big( e_m(\tvr, \tvt) - \tvt s_m(\tvr, \tvt) + \frac{p_m(\tvr, \tvt)}{\tvr} \Big) \vr }, 
\]
and consequently, to obtain the final form of the relative energy inequality: 
\begin{align}
	& \left[ \intO{ E \left(\vr, \vt, \vu \Big| \tvr, \tvt, \tvu \right) } \right]_{t = 0}^{t = \tau} \br 
	&+ \int_0^\tau \intO{ \frac{\tvt}{\vt} \left( \mathbb{S} (\vt, \Ds \vu) : \Ds \vu - \frac{\vc{q}(\vt, \Grad \vt)}{\vt} \right) } \dt \br 
	&\leq  \int_0^\tau \intO{ \frac{\vr}{\tvr} (\vu - \tvu ) \cdot \Grad p(\tvr, \tvt) } \dt \br 
	&- \int_0^\tau \intO{ \left( \vr (s - s(\tvr, \tvt)) \partial_t \tvt + \vr (s - s(\tvr, \tvt)) \vu \cdot \Grad \tvt \right) } \dt  \br &- \int_0^\tau \intO{  
		\left( \frac{\vc{q}(\vt, \Grad \vt)}{\vt} \right) \cdot \Grad \tvt  } \dt \br 
	&- \int_0^\tau \intO{ \Big[ \vr (\vu - \tvu) \otimes (\vu - \tvu) + p(\vr, \vt) \mathbb{I} - \mathbb{S}(\vt, \Ds \vu) \Big] : \Ds \tvu } \dt \br 
	&- \int_0^\tau \intO{ \vr \left[ \partial_t \tvu + (\tvu \cdot \Grad) \tvu + \frac{1}{\tvr} \Grad p(\tvr, \tvt) \right] \cdot (\vu - \tvu) } \dt \br 
	&+ \int_0^\tau \intO{ \left[ \left( 1 - \frac{\vr}{\tvr} \right) \partial_t p(\tvr, \tvt) - \frac{\vr}{\tvr} \vu \cdot \Grad p(\tvr, \tvt) \right] } \dt
	\label{WS8}
\end{align}
for any strong solution $(\tvr, \tvt, \tvu)$ belonging to the regularity 
class \eqref{classS} and emanating from the same initial data as the weak solution $(\vr, \vt, \vu)$.

\begin{Remark} \label{FNR}
The relative energy inequality \eqref{WS8} is the same is its counterpart in 
\cite[Chapter 12, Section 12.3.2]{FeiNovOpen} derived on a bounded domain. Revisiting the proof, we deduce that \eqref{WS8} holds for \emph{any trio} of 
functions $(\tvr, \tvt, \tvu)$ belonging to class \eqref{classS}. 
	
	\end{Remark}

\subsection{Weak--Strong uniqueness, Gronwall's argument}

The proof of Theorem \ref{TM2} will be completed by applying Gronwall's argument 
to \eqref{WS8}. As observed in the preceding part, the inequality 
\eqref{WS8} coincides with its counterpart in \cite[Section 12.3.2]{FeiNovOpen}. 
Consequently, we can repeat step  by step the arguments of \cite[Section 12.3.3]{FeiNovOpen} to arrive at the inequality
\begin{align}
	&\intO{ E \left(\vr, \vt, \vu \Big| \tvr, \tvt, \tvu \right) (\tau, \cdot) }  
	\br 
	&+ \int_0^\tau \intO{ \left( \frac{\tvt}{\vt} - 1 \right) \mathbb{S} (\vt, \Ds \vu) : \Ds \vu  } \dt 
	+ \int_0^\tau \intO{ 
		\left(\frac{\vt}{\tvt} - 1 \right)  \mathbb{S} (\tvt, \Ds \tvu) : \Ds \tvu } \dt \br
	&+ \int_0^\tau \intO{ \Big( \mathbb{S} (\tvt, \Ds \tvu ) - \mathbb{S} (\vt, \Ds \vu) \Big): \Ds (\tvu - \vu) } \dt \br	
	&+\int_0^\tau \intO{ \left( 1 - \frac{\tvt}{\vt} \right) \left( \frac{\vc{q}(\vt, \Grad \vt)\cdot \Grad \vt}{\vt}\right) }  \dt
	\br&	+ \int_0^\tau \intO{ 
		\left(1 - \frac{\vt}{\tvt} \right) \frac{\vc{q} (\tvt, \Grad \tvt) }{\tvt} } \dt
	+ 
	\int_0^\tau \intO{ \left( \frac{ \vc{q}(\vt, \Grad \vt)}{\vt} - \frac{ \vc{q}(\tvt, \Grad \tvt)}{\tvt} \right)  \cdot \Grad (\tvt - \vt)}  \dt \br
	&\leq  
	\int_0^\tau \intO{ F_1 } \dt,
	\label{WS10}
\end{align}	
with a quadratic ``error'' term
\begin{align}
	F_1=  &	\vr (\tvu - \vu) \otimes (\vu - \tvu) : \Ds \tvu 
	+ \left( \frac{\vr}{\tvr} - 1 \right) \Big( \Div \mathbb{S} (\tvt, \Ds \tvu) 
	- \Grad p(\tvr, \tvt \Big) \cdot (\tvu - \vu) \br 
	& + \vr (s(\vr, \vt) - s(\tvr, \tvt)) (\tvu - \vu) \cdot \Grad \tvt \br 
	&- (\vr - \tvr)(s(\vr, \vt) - s(\tvr, \tvt))(\partial_t \tvt + \tvu \cdot \Grad \tvt) + \left(1 - \frac{\vr}{\tvr} \right) (\vu - \tvu) \cdot \Grad p(\tvr, \tvt) \br 
	&+ \Div \tvu \Big( p(\tvr, \tvt)  -\frac{\partial p(\tvr, \tvt)}{\partial \vr} (\tvr -\vr ) -\frac{\partial p(\tvr, \tvt)}{\partial \vt} (\tvt -\vt )- p(\vr, \vt) \Big) \br 
	&+ \tvr \left( s(\tvr, \tvt) - \frac{\partial s(\tvr, \tvt)}{\partial \vr}(\tvr - \vr) - \frac{\partial s(\tvr, \tvt)}{\partial \vt}(\tvt - \vt) - s(\vr, \vt) \right) (\partial_t \tvt + \tvu \cdot \Grad \tvt ).
	\label{WS11}
	\end{align} 
	
Given the convexity of the relative energy, certain error terms in \eqref{WS11} can be immediately absorbed. Accordingly, \eqref{WS10} can be written as 
\begin{align}
	& \intO{ E \left(\vr, \vt, \vu \Big| \tvr, \tvt, \tvu \right) (\tau, \cdot) }  
	\br 
	&+ \int_0^\tau \intO{ \left( \frac{\tvt}{\vt} - 1 \right) \mathbb{S} (\vt, \Ds \vu) : \Ds \vu  } \dt 
	+ \int_0^\tau \intO{ 
		\left(\frac{\vt}{\tvt} - 1 \right)  \mathbb{S} (\tvt, \Ds \tvu) : \Ds \tvu } \dt \br
	&+ \int_0^\tau \intO{ \Big( \mathbb{S} (\tvt, \Ds \tvu ) - \mathbb{S} (\vt, \Ds \vu) \Big): \Ds (\tvu - \vu) } \dt \br	
	&+\int_0^\tau \intO{ \left( 1 - \frac{\tvt}{\vt} \right) \left( \frac{\vc{q}(\vt, \Grad \vt)\cdot \Grad \vt}{\vt}\right) }  \dt
	\br&	+ \int_0^\tau \intO{ 
		\left(1 - \frac{\vt}{\tvt} \right) \frac{\vc{q} (\tvt, \Grad \tvt) }{\tvt} } \dt
	+ 
	\int_0^\tau \intO{ \left( \frac{ \vc{q}(\vt, \Grad \vt)}{\vt} - \frac{ \vc{q}(\tvt, \Grad \tvt)}{\tvt} \right)  \cdot \Grad (\tvt - \vt)}  \dt \br
	&\aleq \int_0^\tau \intO{ E \left(\vr, \vt, \vu \Big| \tvr, \tvt, \tvu \right) } \dt +  
	\left| \int_0^\tau \intO{ F_2 } \dt \right|,
	\label{WS12}
\end{align}	
with 
\begin{align}
	F_2=  &	
	 \left( {\vr} - \tvr \right) \Big( \partial_t \tvu + \tvu \cdot \Grad \tvu 
	\Big) \cdot (\tvu - \vu) \br 
	& + \vr (s(\vr, \vt) - s(\tvr, \tvt)) (\tvu - \vu) \cdot \Grad \tvt \br 
	&- (\vr - \tvr)(s(\vr, \vt) - s(\tvr, \tvt))(\partial_t \tvt + \tvu \cdot \Grad \tvt) + \left(1 - \frac{\vr}{\tvr} \right) (\vu - \tvu) \cdot \Grad p(\tvr, \tvt) \br 
	&+ \Div \tvu \Big( p(\tvr, \tvt)  -\frac{\partial p(\tvr, \tvt)}{\partial \vr} (\tvr -\vr ) -\frac{\partial p(\tvr, \tvt)}{\partial \vt} (\tvt -\vt )- p(\vr, \vt) \Big) \br
	&+ \tvr \left( s(\tvr, \tvt) - \frac{\partial s(\tvr, \tvt)}{\partial \vr}(\tvr - \vr) - \frac{\partial s(\tvr, \tvt)}{\partial \vt}(\tvt - \vt) - s(\vr, \vt) \right) (\partial_t \tvt + \tvu \cdot \Grad \tvt ).
	\label{WS13}
\end{align} 

\subsubsection{Thermal dissipation} 

The following result was proved in \cite[Lemma 4.4]{FeGwKwSG}. 

\begin{Lemma} \label{LWS1}
	Under the hypotheses of Theorem \ref{TM2}, there is a constant $\xi > 0$, 
	\begin{align}
		&\| \vt - \tvt \|_{W^{1,2}(\Omega;\R^3)}^2 \br 	
		&\aleq \intO{ \left( \left( \frac{\tvt}{\vt} - 1 \right)  \frac{\kappa (\vt) |\Grad \vt|^2 }{\vt} + \left(  \frac{\vt}{\tvt} - 1 \right)
			\frac{ \kappa (\tvt) |\Grad \tvt|^2 }{\tvt}  \right)  }  \br
		&+  \intO{ \left(  \frac{\kappa (\tvt) \Grad \tvt }{\tvt} - \frac{\kappa (\vt) \Grad \vt }{\vt}  \right) \cdot
			(\Grad \tvt - \Grad \vt) } + \xi \intO{ E \left( \vr, \vt, \vu \Big| \tvr, \tvt, \tvu \right)	}
		\label{ws18}	
	\end{align}
	where $\xi > 0$ depends on $\sup \tvt$, $\inf \tvt$, $\| \Grad \tvt \|_{L^\infty(\Omega)}$ and $\| [\vt - 2 \vt_\infty]^+ \|_{L^4(\Omega)}$.
	
\end{Lemma}

The result is proved in \cite{FeGwKwSG} on a bounded domain, however, the proof 
on an exterior domain with a positive far field temperature $\vt_\infty$ is the same. 

In view of Lemma \ref{LWS1}, inequality \eqref{WS12} reduces to
\begin{align}
	& \intO{ E \left(\vr, \vt, \vu \Big| \tvr, \tvt, \tvu \right) (\tau, \cdot) }  
	\br 
	&+ \int_0^\tau \intO{ \left( \frac{\tvt}{\vt} - 1 \right) \mathbb{S} (\vt, \Ds \vu) : \Ds \vu  } \dt 
	+ \int_0^\tau \intO{ 
		\left(\frac{\vt}{\tvt} - 1 \right)  \mathbb{S} (\tvt, \Ds \tvu) : \Ds \tvu } \dt \br
	&+ \int_0^\tau \intO{ \Big( \mathbb{S} (\tvt, \Ds \tvu ) - \mathbb{S} (\vt, \Ds \vu) \Big): \Ds (\tvu - \vu) } \dt \br	
	&+ \int_0^\tau \| \vt - \tvt \|^2_{W^{1,2}(\Omega)} \dt \br
	&\aleq \int_0^\tau \intO{ E \left(\vr, \vt, \vu \Big| \tvr, \tvt, \tvu \right) } \dt +  
	\left| \int_0^\tau \intO{ F_2 } \dt \right|.
	\label{WS14}
\end{align}	

\subsubsection{Viscous dissipation}

Now, we may repeat the arguments of \cite[Chapter 4, Section 4.2.3]{FeiNovOpen} to show the inequality
\begin{align}
\| \vu - \tvu \|^2_{D^{1,2}_0(\Omega;\R^3)} &\aleq
\int_0^\tau \intO{ \left( \frac{\tvt}{\vt} - 1 \right) \mathbb{S} (\vt, \Ds \vu) : \Ds \vu  } \dt \br 
&+ \int_0^\tau \intO{ 
	\left(\frac{\vt}{\tvt} - 1 \right)  \mathbb{S} (\tvt, \Ds \tvu) : \Ds \tvu } \dt \br
&+ \int_0^\tau \intO{ \Big( \mathbb{S} (\tvt, \Ds \tvu ) - \mathbb{S} (\vt, \Ds \vu) \Big): \Ds (\tvu - \vu) } \dt \br 
&+ \int_0^\tau \intO{ E \left(\vr, \vt, \vu \Big| \tvr, \tvt, \tvu \right) } \dt.
\label{WS15}
\end{align}	
Similarly to Lemma \ref{LWS1}, the extension of the original proof to an exterior domain is straightforward.

Going back to \eqref{WS14} we conclude 
\begin{align}
	& \intO{ E \left(\vr, \vt, \vu \Big| \tvr, \tvt, \tvu \right) (\tau, \cdot) }  + \int_0^\tau \| \vu - \tvu \|^2_{D^{1,2}_0(\Omega;\R^3)} \dt 
	+ \int_0^\tau \| \vt - \tvt \|^2_{W^{1,2}(\Omega)} \dt \br
	&\aleq \int_0^\tau \intO{ E \left(\vr, \vt, \vu \Big| \tvr, \tvt, \tvu \right) } \dt +  
	\left| \int_0^\tau \intO{ F_2 } \dt \right|.
	\label{WS17}
\end{align}

\subsubsection{Quadratic error terms}

First, we establish some coercivity properties of the relative energy. 
It follows from the structural properties of the function $P$ stated in \eqref{P1}, \eqref{P2}, that the pressure $p_m$ can be written in the form 
\[
p_m(\vr, \vt) = p_{\rm ell}(\vr) + p_M (\vr, \vt), \ \mbox{where}\ p_{\rm ell} = b \vr^{\frac{5}{3}},\ b > 0,  
\]
and 
\[
p_M(\vr, \vt) = \vt^{\frac{5}{2}} P_M(Z),\ Z = \frac{\vr}{\vt^{\frac{3}{2}}}, 
\]
enjoys the same structural properties as $P$, meaning \eqref{c2}--\eqref{c2a}. 

Consequently, evoking the estimates \eqref{P6} we conclude 
\begin{align}
E \left( \vr, \vt, \vu \Big| \tvr, \tvt, \tvu \right) = 
E_M \left( \vr, \vt, \vu \Big| \tvr, \tvt, \tvu \right) + 
H(\vt , \tvt) + \frac{3}{2} \Big( p_{\rm ell} (\vr) - p'_{\rm ell}(\tvr) (\vr - 
\tvr) - p_{\rm ell}(\tvr) \Big), 
	\label{WS18}
	\end{align}
where $H$ enjoys the coercivity properties stated in \eqref{P6}, and $E_M$ is the 
relative energy based on the pressure $p_M$, with the related quantities $e_M$, and $s_M$.

We proceed by estimating the integrals of the error terms $F_2$. 

\bigskip 

\noindent $\bullet$ {\bf Step 1.}
Thanks to the decomposition \eqref{WS18}, we may use the arguments of 
\cite[Section 5.3]{FeiNov2021NON} to control the error integral 
\[
\int_0^\tau \intO{ \left( {\vr} - \tvr \right) \Big( \partial_t \tvu + \tvu \cdot \Grad \tvu 
\Big) \cdot (\tvu - \vu) }. 
\] 
Similarly, we can estimate the integral containing the pressure 
\[
\int_0^\tau \intO{ \Div \tvu \Big( p(\tvr, \tvt)  -\frac{\partial p(\tvr, \tvt)}{\partial \vr} (\tvr -\vr ) -\frac{\partial p(\tvr, \tvt)}{\partial \vt} (\tvt -\vt )- p(\vr, \vt) \Big) } \dt.
\]
Accordingly, inequality \eqref{WS17} reduces to 
\begin{align}
	& \intO{ E \left(\vr, \vt, \vu \Big| \tvr, \tvt, \tvu \right) (\tau, \cdot) }  + \int_0^\tau \| \vu - \tvu \|^2_{D^{1,2}_0(\Omega;\R^3)} \dt 
	+ \int_0^\tau \| \vt - \tvt \|^2_{W^{1,2}(\Omega)} \dt \br
	&\aleq \int_0^\tau \intO{ E \left(\vr, \vt, \vu \Big| \tvr, \tvt, \tvu \right) } \dt +  
	\left| \int_0^\tau \intO{ F_3 } \dt \right|,
	\label{WS19}
\end{align}
with 
\begin{align}
	F_3 &=  	
		 \vr (s(\vr, \vt) - s(\tvr, \tvt)) (\tvu - \vu) \cdot \Grad \tvt \br 
	&- (\vr - \tvr)(s(\vr, \vt) - s(\tvr, \tvt))(\partial_t \tvt + \tvu \cdot \Grad \tvt) + \left(1 - \frac{\vr}{\tvr} \right) (\vu - \tvu) \cdot \Grad p(\tvr, \tvt) \br 
		&+ \tvr \left( s(\tvr, \tvt) - \frac{\partial s(\tvr, \tvt)}{\partial \vr}(\tvr - \vr) - \frac{\partial s(\tvr, \tvt)}{\partial \vt}(\tvt - \vt) - s(\vr, \vt) \right) (\partial_t \tvt + \tvu \cdot \Grad \tvt )\br
  &=F_{3,1}+F_{3,2}+F_{3,3}+F_{3,4}.
	\label{WS20}
\end{align}

\bigskip

\noindent $\bullet$ {\bf Step 2.} We will now estimate the integrals of the error terms in $F_3$. For this we would normally consider different subsets of $\Omega$, depending on the values of $\vr$ and $\vt$. Majority of the cases would be treated similarly as in \cite[Chapter 4, Section 4.2]{FeiNovOpen} .   The new problem solved in the present paper is linked to our choice of the far field boundary conditions, causing that the set where the density is small (relatively to the temperature, i.e. $Z$ is small) might be unbounded. 

In accordance with this, from now on we will focus solely on that case, restricting the analysis to the part of the physical space, where 
\[
Z = \frac{\vr}{\vt^{\frac{3}{2}}} \leq D, 
\]
for a large $D$. Thus we can fix $D$ so that 
\begin{equation} \label{WS24}
\frac{\vr}{\vt^{\frac{3}{2}}}, \ \frac{\tvr}{\tvt^{\frac{3}{2}}} \leq D.
\end{equation}

Consequently, there exists $\Lambda > 0$ such that 
\[
P_M(Z) = \Lambda Z + P_B(Z) ,
\] 
where $P_B(Z)$ complies with the thermodynamic stability hypothesis formulated 
in \eqref{c2} for $0 \leq Z \leq D$. Thus we may write the relative energy \eqref{WS18} once more in the form 
\begin{align}
	E \left( \vr, \vt, \vu \Big| \tvr, \tvt, \tvu \right) &= 
	E_B \left( \vr, \vt, \vu \Big| \tvr, \tvt, \tvu \right) + 
	H(\vt , \tvt) + \frac{3}{2} \Big( p_{\rm ell} (\vr) - p'_{\rm ell}(\tvr) (\vr - 	\tvr) - p_{\rm ell}(\tvr) \Big) \br 
	&+ E_\Lambda \left( \vr, \vt, \vu \Big| \tvr, \tvt, \tvu \right), 
	\label{WS25}
\end{align}
where $E_\Lambda$ is associated to the linear function $P(Z) = \Lambda Z$ as long as the parameters $(\vr, \tvr, \vt, \tvt)$ belong to the range \eqref{WS24}. 

The advantage of \eqref{WS25} is that the function $E_\Lambda$ can be explicitly computed. Specifically, 
\[
p = \Lambda \vr \vt,\ e = \Lambda \frac{3}{2} \vt,\ 
s = \Lambda \left( \frac{3}{2} \log (\vt) - \log(\vr) \right).
\]
Accordingly, the associated relative energy reads
\begin{align}
E_\Lambda &\left( \vr, \vt, \vu \Big| \tvr, \tvt, \tvu \right) \br &= 
\Lambda \left[ \frac{3}{2} \vr \vt - \tvt \vr  \left( \frac{3}{2} \log (\vt) - \log(\vr) \right) - \vr \left( \frac{5}{2} \tvt - \tvt \left( \frac{3}{2} \log (\tvt) - \log(\tvr) \right)    \right)   + \tvr \tvt \right] \br 
&= \frac{3}{2} \Lambda \vr \left[(\vt - \tvt) - \tvt \Big( \log(\vt) - \log(\tvt) \Big)      \right]  + \Lambda \tvt \left[  \vr \Big( \log(\vr)  -   \log(\tvr) \Big)   + \tvr  - \vr  \right]. 
\label{WS26}
\end{align}
Note that both 
\[
\vt \mapsto \left[(\vt - \tvt) - \tvt \Big( \log(\vt) - \log(\tvt) \Big)      \right]
\]
and
\[
\vr \mapsto \left[  \vr \Big( \log(\vr)  -   \log(\tvr) \Big)   + \tvr  - \vr  \right]
\]
are strictly convex functions attaining their minimum at $\tvt$ and $\tvr$, respectively.

As the function $\tvt$ is bounded below away from zero, we get 
\begin{align}
\Lambda \tvt \left[  \vr \Big( \log(\vr)  -   \log(\tvr) \Big)   + \tvr  - \vr  \right] &\ageq \frac{1}{\vr} (\vr - \tvr)^2 \ \mbox{if}\ \vr \geq \tvr,  \br 
\Lambda \tvt \left[  \vr \Big( \log(\vr)  -   \log(\tvr) \Big)   + \tvr  - \vr  \right] &\ageq \frac{1}{\tvr} (\vr - \tvr)^2 \ \mbox{if}\ \vr \leq \tvr.
 \label{WS27}
\end{align}
\bigskip
The following steps are devoted to control of integrals of terms $F_{3,1}, F_{3,2}, F_{3,3}, F_{3,4}$ from \eqref{WS20}. 
In agreement with our previous discussion, it is enough to consider the 
``essential set''
\begin{equation} \label{WS28}
\mathcal{O}_{\rm ess} = \left\{ 	(t,x) \in (0,T) \times \Omega \ \Big| \ \ 0 \leq \vr(t,x) \leq \Ov{\vr},\ 
	0 < \underline{\vt} \leq \vt(t,x) \leq \Ov{\vt} \right\}, 
	\end{equation}
where	
\[
\Ov{\vr} > 2 \sup_{(0,T) \times \Omega} \tvr,\  \quad 
\underline{\vt} \leq \frac{1}{2} \inf_{(0,T) \times \Omega} \tvt ,\quad \ 
\Ov{\vt} \geq 2 \sup_{(0,T) \times \Omega} \tvt. 
\]

\noindent $\bullet$ {\bf Step 3: Estimate of $F_{3,1}$ and of of $F_{3,3}$.} 

As for the integral of the first term in \eqref{WS20}
\[
\int_0^\tau \intO{ \mathds{1}_{\rm \mathcal{O}_{\rm ess}} F_{3,1}}=\int_0^\tau \intO{ \mathds{1}_{\rm \mathcal{O}_{\rm ess}} \vr \Big( s(\vr, \vt) - s(\tvr, \tvt) \Big) (\tvu - \vu) \cdot \Grad \tvt } \dt, 
\]  
we first handle the radiation component 
\begin{align}
 &\intO{ \CH \vr \Big( \frac{\vt^3}{\vr} - \frac{\tvt^3}{\tvr} \Big) (\tvu - \vu) \cdot \Grad \tvt }  = 
 \intO{ \CH \sqrt{\tvr}\Big( \vt^3 - \tvt^3 \Big) \sqrt{\tvr}(\tvu - \vu) \cdot \frac{\Grad \tvt}{\tvr} }  \br
&+  \intO{ \CH \tvt^3 \Big( \vr - \tvr \Big) (\tvu - \vu) \cdot \frac{\Grad \tvt}{\tvr} }. 
\label{WS28b}
\end{align}
In view of  \eqref{proper}, we have $
\left| \frac{\Grad \tvt}{\tvr} \right| \aleq 1 .$
Consequently, in accordance with \eqref{proper} and by the
Hardy--Sobolev inequality, 
\begin{align}
&\intO{ \CH \sqrt{\tvr}\Big( \vt^3 - \tvt^3 \Big) \sqrt{\tvr}(\tvu - \vu) \cdot \frac{\Grad \tvt}{\tvr} } \leq 
\delta \intO{ \tvr |\tvu - \vu|^2 } + 
C(\delta) \| \vt - \tvt \|^2_{L^2(\Omega)} \br 
&\quad \leq \delta \intO{ \frac{ |\tvu - \vu|^2 }{|x|^2} } + 
c(\delta) \intO{ E \left(\vr, \vt, \vu \Big| \tvr, \tvt, \tvu \right) } \br
&\quad \leq \delta \| \tvu - \vu \|^2_{D^{1,2}_0(\Omega; R^3)} + 
c(\delta) \intO{ E \left(\vr, \vt, \vu \Big| \tvr, \tvt, \tvu \right) }.
\label{WS28c}
\end{align}
As for the second integral in \eqref{WS28b}, we split it in two parts,
\begin{align}
\intO{ \CH \tvt^3 \Big( \vr - \tvr \Big) (\tvu - \vu) \cdot \frac{\Grad \tvt}{\tvr} } &\aleq \intO{ \CH \mathds{1}_{ \{ \vr \geq \tvr \} } | \vr - \tvr | \ | \tvu - \vu|} \br & +\intO{ \CH \mathds{1}_{ \{ \vr \leq \tvr \} } | \vr - \tvr | \ | \tvu - \vu|} .
	\nonumber
	\end{align}
In accordance with \eqref{WS27}, the first part is estimated as follows	
\begin{align}
\intO{ \CH \mathds{1}_{ \{ \vr \geq \tvr \} } | \vr - \tvr | \ | \tvu - \vu|} 
&\aleq \intO{ \CH \mathds{1}_{  \{ \vr \geq \tvr \} } \frac{|\vr - \tvr|^2 }{\vr} } 
+ \intO{ \vr |\vu - \tvu|^2 } \br 
&\aleq \intO{ E \left(\vr, \vt, \vu \Big| \tvr, \tvt, \tvu \right) },
\nonumber
\end{align}
while for the second part, similarly as in \eqref{WS28c} we write
\begin{align}
\intO{ \CH \mathds{1}_{ \{ \vr \leq \tvr \} } | \vr - \tvr | \ | \tvu - \vu|} 
&\aleq C(\delta)\intO{ \CH \mathds{1}_{  \{ \vr \leq \tvr \} } \frac{|\vr - \tvr|^2 }{\tvr} } 
+ \delta\intO{ \tvr |\vu - \tvu|^2 } \br 
&\aleq C(\delta) \intO{ E \left(\vr, \vt, \vu \Big| \tvr, \tvt, \tvu \right) }+\delta \| \tvu - \vu \|^2_{D^{1,2}_0(\Omega; R^3)}.
\nonumber
\end{align}

Note in passing that the integral 
\[ \int_0^\tau \intO{  \CH F_{3,3}}=
\int_0^\tau \intO{ \CH \left(1 - \frac{\vr}{\tvr} \right) (\vu - \tvu) \cdot \Grad p(\tvr, \tvt) } \dt
\]
can be handled in the same way. 

Thus we are left with the non-radiation component of $F_{3,1}$
\[
\int_0^\tau \intO{ \mathds{1}_{\rm \mathcal{O}_{\rm ess}} \vr \Big( \mathcal{S}(Z) - \mathcal{S}(\widetilde{Z}) \Big) (\tvu - \vu) \cdot \Grad \tvt } \dt,\quad \text{where}\  Z = \frac{\vr}{\vt^{\frac{3}{2}}},\ \widetilde{Z} = \frac{\tvr}{\tvt^{\frac{3}{2}}}.
\]
As $\mathcal{S}'(Z) \approx \frac{1}{Z}$ for $Z \aleq 1$, we have 
\begin{equation} \label{WS29}
\left| \mathcal{S}(Z) - \mathcal{S}(\widetilde{Z}) \right| \aleq \left| 
\log(Z) - \log (\widetilde{Z}) \right|. 
\end{equation}
Consequently, 
\begin{align} 
	&\left| \int_0^\tau \intO{ \mathds{1}_{\rm \mathcal{O}_{\rm ess}} \vr \Big( \mathcal{S}(Z) - \mathcal{S}(\widetilde{Z}) \Big) (\tvu - \vu) \cdot \Grad \tvt } \dt \right| \br &\quad \aleq 
	\int_0^\tau \intO{ \mathds{1}_{\rm \mathcal{O}_{\rm ess}} \sqrt{\vr} \Big| \log(\vt) - \log(\tvt) \Big| \sqrt{\vr} |\tvu - \vu||\Grad \tvt| } \dt \br 
	&\quad + \int_0^\tau \intO{ \mathds{1}_{\rm \mathcal{O}_{\rm ess}} \vr \Big| \log(\vr) - \log(\tvr) \Big|  |\tvu - \vu||\Grad \tvt| } \dt, 
	\nonumber
	\end{align}
where, obviously, 
\[
\int_0^\tau \intO{ \mathds{1}_{\rm \mathcal{O}_{\rm ess}} \sqrt{\vr} \Big| \log(\vt) - \log(\tvt) \Big| \sqrt{\vr} |\tvu - \vu||\Grad \tvt| } \dt \aleq 
\int_0^\tau \intO{ E \left( \vr, \vt, \vu \Big| \tvr, \tvt, \tvu \right) } \dt.
\]

Finally, 
\begin{align}
&\int_0^\tau \intO{ \mathds{1}_{\rm \mathcal{O}_{\rm ess}} \vr \Big| \log(\vr) - \log(\tvr) \Big|  |\tvu - \vu||\Grad \tvt| } \dt \br 
&\quad \aleq \int_0^\tau \intO{ \mathds{1}_{\rm \mathcal{O}_{\rm ess}} \tvr \Big[ \vr ( \log(\vr) - \log(\tvr))  + \tvr - \vr \Big]  |\tvu - \vu| \frac{|\Grad \tvt|}{\tvr} } \dt \br 
&\quad + \int_0^\tau \intO{ |\vr - \tvr| |\vu - \tvu| |\Grad \tvt| } \dt.
\nonumber
\end{align}
By virtue of the H\" older and the Sobolev inequalities, we get
\begin{align}
	&\intO{ \mathds{1}_{\rm \mathcal{O}_{\rm ess}} \tvr \Big[ \vr ( \log(\vr) - \log(\tvr))  + \tvr - \vr \Big]  |\tvu - \vu| \frac{|\Grad \tvt|}{\tvr} }  \br 
	&\quad \leq \delta \| \vu - \tvu \|^2_{L^6(\Omega; R^3)}  + C(\delta) \left(  \intO{ \CH \tvr^{\frac{6}{5}  }\Big[ \vr ( \log(\vr) - \log(\tvr))  + \tvr - \vr \Big]^{\frac{6}{5}}}  \right)^{\frac{5}{3}}.   
	\nonumber
\end{align}
Moreover, 
\begin{align} 
&\left(  \intO{ \CH \tvr^{\frac{6}{5}}  \Big[ \vr ( \log(\vr) - \log(\tvr))  + \tvr - \vr \Big]^{\frac{6}{5}} } \right)^{\frac{5}{3}} \br
&\quad \aleq \left(  \intO{  \tvr^{\frac{6}{5}}  \Big[ \vr ( \log(\vr) - \log(\tvr))  + \tvr - \vr \Big]^{\frac{6}{5}- 1} \Big[ \vr ( \log(\vr) - \log(\tvr))  + \tvr - \vr \Big] }\right)^{\frac{5}{3}} \br 
&\quad \aleq  \left(  \intO{ \CH  \Big[ \vr ( \log(\vr) - \log(\tvr))  + \tvr - \vr \Big]} \right)^{\frac{5}{3}} \aleq 
\left( \intO{ E_\Lambda \left(\vr, \vt, \tvt \Big| \tvr, \tvt.  \tvu \right) } \right)^{\frac{5}{3}} \br 
&\quad \aleq 
 \intO{ E_\Lambda \left(\vr, \vt, \tvt \Big| \tvr, \tvt, \tvu \right) }.
\nonumber 	
\end{align}

\bigskip

\noindent $\bullet$ {\bf Step 4: Estimate of $F_{3,2}$ and $F_{3,4}$.} Finally, we have to control the entropy dependent terms in \eqref{WS20}, namely 	
\begin{align}
(\tvr - \vr)&(s(\vr, \vt) - s(\tvr, \tvt))(\partial_t \tvt + \tvu \cdot \Grad \tvt)
\br 
&+ \tvr \left( s(\tvr, \tvt) - \frac{\partial s(\tvr, \tvt)}{\partial \vr}(\tvr - \vr) - \frac{\partial s(\tvr, \tvt)}{\partial \vt}(\tvt - \vt) - s(\vr, \vt) \right) (\partial_t \tvt + \tvu \cdot \Grad \tvt ) \br 
&= \left[ \vr (s(\tvr, \tvt) - s(\vr, \vt))  - \tvr \left(           \frac{\partial s(\tvr, \tvt)}{\partial \vr}(\tvr - \vr) + \frac{\partial s(\tvr, \tvt)}{\partial \vt}(\tvt - \vt)    \right)           \right] \Big(\partial_t \tvt + \tvu \cdot \Grad \tvt \Big).
\nonumber
\end{align}

Consider the function 
\[
\mathcal{F}(\vr, \vt) \mapsto \left[ \vr (s(\tvr, \tvt) - s(\vr, \vt))  - \tvr \left(           \frac{\partial s(\tvr, \tvt)}{\partial \vr}(\tvr - \vr) + \frac{\partial s(\tvr, \tvt)}{\partial \vt}(\tvt - \vt)    \right)           \right].
\]
A direct computation yields 
\[
\mathcal{F}(\tvr, \tvt) = 0, 
\]
\[
\frac{\partial}{\partial \vr} \mathcal{F}(\vr, \vt) = 
 s(\tvr, \tvt) - s(\vr, \vt) - \vr \frac{\partial }{\partial \vr} s(\vr, \vt) + 
\tvr \frac{\partial }{\partial \vr} s(\tvr, \tvt) \ \Rightarrow\ 
\frac{\partial}{\partial \vr} \mathcal{F}(\tvr, \tvt) = 0,
\]
and
\[
\frac{\partial}{\partial \vt} \mathcal{F}(\vr, \vt) =  
- \vr \frac{\partial }{\partial \vt} s(\vr, \vt) + \tvr \frac{\partial }{\partial \vt} s(\tvr, \tvt) \ \Rightarrow\ 
\frac{\partial}{\partial \vt} \mathcal{F}(\tvr, \tvt) = 0.
\]
In particular, we can discard the radiation component of the entropy as 
the corresponding part of $\mathcal{F}$ depends only on $\vt$ and is therefore regular. Accordingly, 
we focus on 
\[
\mathcal{F}_m(\vr, \vt) \mapsto \left[ \vr (s_m(\tvr, \tvt) - s_m(\vr, \vt))  - \tvr \left(           \frac{\partial s_m(\tvr, \tvt)}{\partial \vr}(\tvr - \vr) + \frac{\partial s_m(\tvr, \tvt)}{\partial \vt}(\tvt - \vt)    \right)           \right].
\]

Similarly to the above, a direct computation yields
\begin{align} 
\frac{\partial^2 }{\partial \vr^2} \mathcal{F}_m (\vr, \vt) &= 
- 2 \frac{\partial }{\partial \vr} s_m(\vr, \vt) - \vr \frac{\partial^2 }{\partial \vr^2} s_m(\vr, \vt), \br 
\frac{\partial^2 }{\partial \vr \partial \vt} \mathcal{F}_m (\vr, \vt) &=
- \frac{\partial }{\partial \vt} s_m(\vr, \vt) - \vr 
\frac{\partial^2 }{\partial \vr \partial \vt} s_m(\vr, \vt), \br 
\frac{\partial^2 }{\partial \vt^2} \mathcal{F}_m (\vr, \vt) &= - 
\vr \frac{\partial^2 }{\partial \vt^2} s_m (\vr, \vt).
\nonumber
	\end{align}
	
Now, in accordance with hypothesis \eqref{c2b}, $P''(0) = 0$. In particular, the 
entropy $s_m$ can be calculated for the \emph{linear} function $P(Z) = P'(0) Z=\Lambda Z $, 
\[
s_m(\vr, \vt) = \Lambda \left( \frac{3}{2} \log(\vt) - \log(\vr) \right) + 
s_b (\vr, \vt), 
\]
where $s_b \in C^2[0, \infty)$, whence its contribution to $\mathcal{F}$ can be controlled. The assuming $s_m \approx \left( \frac{3}{2} \log(\vt) - \log(\vr) \right)$ we get 
\begin{align} 
	\frac{\partial^2 }{\partial \vr^2} \mathcal{F}_m (\vr, \vt) &= 
	\frac{1}{\vr}, \br 
	\frac{\partial^2 }{\partial \vr \partial \vt} \mathcal{F}_m (\vr, \vt) &=
	\frac{3}{2} \log(\vt) , \br 
	\frac{\partial^2 }{\partial \vt^2} \mathcal{F}_m (\vr, \vt) &=  
	\frac{9}{4} \frac{\vr}{\vt^2} .
	\nonumber
\end{align}
As $\vt$ is bounded below on the essential set, the integral 
\[
\int_0^\tau \intO{\CH\left[ \vr (s(\tvr, \tvt) - s(\vr, \vt))  - \tvr \left(           \frac{\partial s(\tvr, \tvt)}{\partial \vr}(\tvr - \vr) + \frac{\partial s(\tvr, \tvt)}{\partial \vt}(\tvt - \vt)    \right)           \right] \Big(\partial_t \tvt + \tvu \cdot \Grad \tvt \Big)} \dt  
\]
can be controlled by 
\[
\int_0^\tau \intO{ E_\Lambda \left(\vr, \vt, \vu \Big| \tvr, \tvt, \tvu \right) } \dt
\]

We have proved Theorem \ref{TM2}.

\def\cprime{$'$} \def\ocirc#1{\ifmmode\setbox0=\hbox{$#1$}\dimen0=\ht0
	\advance\dimen0 by1pt\rlap{\hbox to\wd0{\hss\raise\dimen0
			\hbox{\hskip.2em$\scriptscriptstyle\circ$}\hss}}#1\else {\accent"17 #1}\fi}

\end{document}